\documentclass[twocolumn]{autart}    

\usepackage{graphicx,color}
\usepackage[dvips]{epsfig}
\usepackage{subfigure}
\usepackage[comma,numbers,square,sort&compress]{natbib}
\usepackage{exscale,amsmath}
\usepackage{amssymb,cases,bigints,setspace}
\usepackage{amsfonts}
\usepackage{epstopdf}
\usepackage{algorithmic}
\usepackage{dsfont}
\usepackage{anyfontsize}

\newtheorem{Theorem}{Theorem}[section]
\newtheorem{Lemma}{Lemma}[section]
\newtheorem{Proposition}{Proposition}[section]
\newtheorem{Remark}{Remark}[section]

\newtheorem{Problem}{Problem}[section]
\begin{document}

\begin{frontmatter}

\title{Linear quadratic mean-field game-team analysis: a mixed coalition approach}

\thanks[footnoteinfo]{Corresponding author Shujun Wang. This work was supported by the National Key R\&D Program of China (No. 2022YFA1006104), RGC P0037987, the National Natural Science Foundation of China (Nos. 12001320, 11831010, 61961160732), the Natural Science Foundation of Shandong Province (No. ZR2019ZD42), the Taishan Scholars Climbing Program of Shandong (No. TSPD20210302), the Taishan Scholars Young Program of Shandong (No. TSQN202211032) and the Young Scholars Program of Shandong University.}

\author[address1]{Jianhui Huang}\ead{majhuang@polyu.edu.hk},
\author[address1]{Zhenghong Qiu}\ead{zhenghong.qiu@connect.polyu.hk},
\author[address2]{Shujun Wang}\ead{wangshujun@sdu.edu.cn},
\author[address3]{Zhen Wu}\ead{wuzhen@sdu.edu.cn}

\address[address1]{Department of Applied Mathematics, The Hong Kong Polytechnic University, Hong Kong}
\address[address2]{School of Management, Shandong University, Jinan, Shandong 250100, China}
\address[address3]{School of Mathematics, Shandong University, Jinan, Shandong 250100, China}


\begin{abstract}
Mean-field theory has been extensively explored in decision analysis of {large-scale} (LS) systems but traditionally in ``pure" cooperative or competitive settings. This leads to the so-called mean-field game (MG) or mean-field team (MT). This paper introduces a new class of LS systems with cooperative inner layer and competitive outer layer, so a ``mixed" mean-field analysis is proposed for distributed game-team strategy. A novel asymptotic mixed-equilibrium-optima is also proposed and verified.
\end{abstract}

\end{frontmatter}

\def\d{\delta}
\def\m{\bar{\mathbf{m}}}

\def\D{\Delta}

\def\P{\mathbb{P}}

\def\mD{\mathcal{D}}

\def\E{\mathbb{E}}

\def\S{\mathbb{S}}

\def\I{\mathcal{I}}

\def\A{\mathcal{A}}

\def\B{\mathcal{B}}
\def\bH{\mathbf{H}}

\def\C{\mathcal{C}}

\def\U{\mathcal{U}}

\def\J{\mathcal{J}}

\def\M{\mathbf{M}}

\def\L{\Lambda}

\def\bL{\mathbf{\Lambda}}

\def\F{\mathcal{F}}

\def\H{\mathcal{H}}

\def\bF{\mathbb{F}}

\def\R{\mathbb{R}}

\def\xN{x^{(N)}}

\def\uN{u^{(N)}}

\def\bmfu{\bar{\mathbf{u}}}

\def\bP{\mathbf{P}}

\def\bA{\mathbf{A}}
\def\bbA{\mathbb{A}}

\def\bB{\mathbf{B}}
\def\bbB{\mathbb{B}}

\def\bC{\mathbf{C}}
\def\bbC{\mathbb{C}}

\def\by{\mathbf{y}}

\def\bD{\mathbf{D}}

\def\bQ{\mathbf{Q}}
\def\bbQ{\mathbb{Q}}

\def\bR{\mathbf{R}}
\def\bbR{\mathbb{R}}

\def\bG{\mathbf{G}}

\def\bx{\mathbf{x}}

\def\bu{\mathbf{u}}

\def\bz{\mathbf{z}}

\def\bV{\mathbf{V}}
\def\bK{\mathbf{K}}

\def\bv{\mathbf{v}}

\definecolor{darkred}{rgb}{0.7,0,0}

\def\r{\color{red}}
\def\dr{\color{darkred}}
\def\b{\color{blue}}
\def\bk{\color{black}}

\def\+{\!+\!}
\def\-{\!-\!}
\def\={\!=\!}
\def\x{\bar{x}}
\def\cx{\check{x}}
\def\tx{\tilde{x}}

\def\u{\bar{u}}
\def\cu{\check{u}}
\def\tu{\tilde{u}}

\def\p{\bar{p}}

\def\spn{\mathrm{span}}

\def\diag{\mathrm{diag}}

\def\s{\small}
\def\f{\footnotesize}
\def\n{\normalsize}

\def\T{\int_{0}^{T}}

\def\Eo{\mathbb{E}_0}

\def\hmc{{\text{\rm\textbf{HMC}}}}

\def\hmg{{\text{\rm\textbf{HMG}}}}

\def\hmt{{\text{\rm\textbf{HMT}}}}

\def\ccmg{{\text{\rm\textbf{CC}$_{\text{\rm\textbf{\tiny MG}}}$}}}

\def\ccmt{{\text{\rm\textbf{CC}$_{\tiny\text{\rm\textbf{\tiny MT}}}$}}}

\def\sumI1{\sum_{k\in\I_1, k\neq i}}

\def\scl8{\scaleobj{.8}}

\section{Introduction}

This paper is inspired by (dynamic) decision of large-scale (LS) systems for which the most salient feature, in linear-quadratic (LQ) setting, is existence of sufficiently many negligible agents interacted among their states or objectives via empirical state-average or control-average. Weakly-coupled LS systems have been found broad applications across economics, finance, biology and engineering. Interested readers may refer \cite{hcm2007a}, \cite{ll2007}, etc. For LS systems with highly complex interactions, mean-field theory provides an effective scheme to study its decision \emph{asymptotically} as population size $N \rightarrow \infty$.

In principle, decisions of LS systems can be classified into \emph{non-cooperative} game or \emph{cooperative} team, relying on coalition structure formalized by involved agents in underlying population. Accordingly, mean-field game (MG) or mean-field team (MT) arise naturally when studying LS asymptotically. Along game-line (MG), all agents, $\{\mathcal{A}_i\}_{i=1}^{N}$, are non-cooperative with competitive objectives $\{\mathcal{J}_{i}\}_{i=1}^{N}$ thus some Nash type equilibrium should be analyzed. There have accumulated vast literature on {MG} study, see e.g. \cite{bsyy2016}, \cite{blp2014}, \cite{c2010}, \cite{cd2013}, \cite{gs2014}, \cite{HHN18}, \cite{hcm2007a}, \cite{ll2007}, \cite{mb2018}, \cite{NNY2020}, \cite{tzb2014}. Parallel to game, MT forms another appealing line along which $\{\mathcal{A}_i\}_{i=1}^{N}$ are cooperative to same team (social) objective $\mathcal{J}_{soc}^{(N)}=\sum \mathcal{J}_{i}$. MT has also been extensively explored from a wide range of perspectives, e.g., \cite{hcm2012}, \cite{ncmh2013}, \cite{nm2018}, \cite{wz2017}.
All aforementioned (MG, MT) attempts, albeit well explored, have only been built on a ``pure" basis whereby all LS agents are purely cooperative or competitive. However, in reality, LS system often displays some ``mixed" behaviors in its agents' organization with combined game and team. To be explained later, pure basis seems over idealized to fit reality, so some expanded analysis on more realistic mixed basis is strongly suggested.

On one hand, from a ``practical" viewpoint, various real environments suggest LS to display some ``mix" decision patterns. Indeed, mixed LS come often from economics, engineering, or management, etc., when two or more competitive networks co-exist, e.g., duopoly market with distributed franchisees; adversarial networks with operation knots (e.g., \cite{BC98}, \cite{DHM20}). Specifically, a mixed basis is well-posed whenever classical \emph{two-person game} or general multiple-person-game are framed with \emph{decentralized} decisions on their distributed \emph{sub-units}. By this, original two (multiple) persons or entities still play game, while all sub-units within each entity formalize two (multiple) teams. So a ``mixed" basis arises with tangled game and team. On the other hand, from a ``mathematical" viewpoint: conceptually, any LS decision can be captured by its underlying \textbf{coalition matrix} (to be introduced in Section \ref{coa}), an $N\times N$ matrix satisfying given block structure properties. Exactly, MG and MT are captured by its two extreme cases: \emph{identical} and \emph{one} matrix, respectively ($\mathbf{C}_1$ and $\mathbf{C}_2$). However, besides them, there exists other types coalition matrix which should feature other meaningful LS decision patterns. More general mixed structure can be identified using other variants of coalition matrix.

Motivated by above discussions, this paper formulates and analyzes a certain type of \emph{mixed} game-team problem involving two levels of interactions: cooperative inner layer for weakly-coupled agents and competitive outer layer between sub-systems. This yields an adversarial LS networks with complex interactions inside and outside. Accordingly, our main contributions are as follows:
\begin{enumerate}
  \item A new \emph{Mix} game-team is introduced and studied, which offers a more general LS setting including those classical ones on mean-field game and team.
  \item A coalition matrix representation is formulated to characterize competitive and cooperative structure among weakly-coupled decision makers. Specifically, classical MG and MT correspond to two extreme cases of such coalition structures.
  \item  A bilateral person-by-person optimality condition and bilateral variation decomposition are proposed for the first time to design distributed strategy. Wellposedness of associated consistency condition (CC) is also discussed under mild conditions.
  \item A new notion of asymptotic mixed-equilibrium-optima is proposed to mixed LS study, and verified to the derived distributed mixed-game-team strategy asymptotically.
\end{enumerate}The remainder of this paper is organized as follows. Section \ref{pre} gives some preliminary notations and a general formulation of mixed game-team via coalition matrix. Section \ref{sec 3} makes some synthesis analysis on variation decomposition. Section \ref{sec 4} aims to some distributed mixed strategy design via a bilateral auxiliary problem and resulting CC system. Asymptotic mixed-equilibrium-optima is studied in Section \ref{sec 5}. Section \ref{conc} concludes this work.

\section{Problem formulation}\label{pre}
\subsection{Preliminary}
Throughout this paper, $\prod$ denotes the Cartesian product. $\mathbb{R}^{n\times m}$ and
$\mathbb{S}^n$ denote the sets of all $(n\times m)$ real matrices and
 all $(n\times n)$ symmetric matrices respectively. For a vector or matrix $A$, $A^T$ denotes transpose of $A$. 
 We denote $\|\cdot\|$ as the standard Euclidean
norm and $\langle\cdot,\cdot\rangle$ as the standard Euclidean inner
product. For a vector $v$ and a symmetric matrix $S$, $\|v\|_S^2:
=\langle Sv,v\rangle\=v^{T}Sv$.  $S> 0 \ (\geq 0)$ means that
$S$ is positive (semi-positive) definite, and $S\gg 0$ means that, $S
- \varepsilon I \geq 0$, for some $\varepsilon>0$.

Moreover, we suppose that $(\Omega,\F,\mathbb{P})$ is an
complete probability space, and $W$ = $(W_{1},\cdots,W_N)$ is a $N$-dimensional standard Brownian
motion defined on it. $\{\F(t)\}_{t\geq 0}$ is the natural filtration generated by $W_i(s), 0\leq s\leq t,1\leq i\leq N$, augmented by all the $\mathbb{P}-$null sets of $\F$, and $\bF:=\{\F(t)\}_{t\geq0}$. $\{\F_i(t)\}_{t\geq 0}$ is the natural filtration generated by $W_i(s), 0\leq s\leq t$, augmented by all the $\mathbb{P}-$null sets of $\F_i$, and $\bF_i:=\{\F_i(t)\}_{t\geq0}$, $1\leq i\leq N$. For a generic Euclidean space $\mathcal{E}$ and filtration $\mathbb{K}$, we introduce the following spaces: $L^{\infty}(0,T;\mathcal{E})=\left\{x : [0,T]\rightarrow \mathcal{E}\big| \ x  \text{ is bounded and}\ \text{deterministic}\right\},$ \\ $C([0,T];\mathcal{E})=\left\{x : [0,T]\rightarrow \mathcal{E}\big| \ x  \text{ is continuous}\right\}$,\\ $L_{\mathbb{K}}^2(0,T;\mathcal{E})=\left\{x : [0,T]\times\Omega\rightarrow \mathcal{E}\big| \ x  \text{ is }\mathbb{K}\text{-progressively}\right.$\\  $\text{measurable, } \|x \|_{L^2} := \Big(\E  {\T}\|x(t)\|^2dt\Big)^{\frac{1}{2}}<\infty\Big\}.$


\subsection{Large-scale system with weak-coupling} \label{moti}
In this paper, on $(\Omega,\F,\bF,\mathbb{P})$ we consider an LS system with $N$ agents, denoted by $\{\A_{i}\}_{i\in\I}$, and $\I=\{1,\cdots,N\}$ denotes the index set of the agents. The aggregation of all  agents is denoted by $\A:=\{\A_{i}\}_{i\in\I}$. The state process of the $i^{\text{th}}$ agent $ \A_i$ is modeled by a controlled linear SDE on finite time horizon $[0,T]$:
\begin{equation}\label{ls} \begin{aligned}
 \mathcal{S}_i:\quad   &dx_i(t)\= (A_ix_i(t) \+ B_iu_i(t) \+ F_ix^{(N)}(t))dt \\  &\qquad\qquad\+  \sigma_i dW_{i}(t),\quad x_i(0)\=\xi_{i},\\
    \end{aligned} \end{equation}
     and for the sake of notation simplicity, we denote  $\bu := (u_1,\cdots,u_N)$ and $\bx := (x_1,\cdots,x_N)$. Then $\mathcal{S} := \{\mathcal{S}_i\}_{i\in\I}$ form a weakly coupled LS  system, since each agent $\A_i$ is coupled with the others via $\xN$, and  $\A_i$ could only provide weak influence to the others  at $O\left(\frac{1}{N}\right)$ order.
To evaluate each control law $u_i$, we introduce the following individual cost functional:
\begin{equation}\label{JMG}
\begin{aligned}
&\J_i(\xi_0;\bu)\=\J_i(\xi_0;{u}_i,{\bu}_{-i})\\
\=&\frac{1}{2}\E \Bigg\{{\T}\Big[\|x_i(t)\-\Gamma_i\xN(t)
\|_{Q_i}^{2} \+ \|u_i(t)\|^{2}_{R_i}\Big] dt \Bigg\},\\
\end{aligned}
\end{equation}
where  $\bu_{-i} = (u_1,\cdots,u_{i-1},u_{i+1},\cdots,  u_N)$.
LS system \eqref{ls}-\eqref{JMG} is in linear quadratic (LQ) setting that is well documented in mean-field literature (e.g. \cite{bsyy2016}, \cite{hcm2007a}, \cite{wz2017}, etc). Henceforward, for notation simplicity, we will drop $t$ if no confusion.
For $i \in \I$, the decentralized admissible strategy set for the $i^{th}$ agent is given by
$\U_i = \Big\{u_i|u_i\in L_{\mathbb{F}_i}^2(0,T;\R^m)\Big\}$; the centralized admissible strategy set for the $i^{th}$ agent is given by
$\U_i^c = \Big\{u_i|u_i\in L_{\mathbb{F}}^2(0,T;\R^m)\Big\}$.

For LS \eqref{ls}, centralized decision on full information $\mathbb{F}$ becomes inefficient due to coupling. Instead, we prefer decentralized one on distributed $\mathbb{F}_{i}$ via mean-field scheme. Along this line, recall the well-studied MG and MT:
\begin{equation}\label{mgp}
\left\{\begin{aligned}
&\text{\textbf{MG:}}\inf_{u_{i}(\cdot) \in \U_{i}}\mathcal{J}_{i}(u_{i}; \textbf{u}_{-i});\\
&\text{\textbf{MT:}}\inf_{\textbf{u}(\cdot) \in \U^N}\mathcal{J}_{soc}^{(N)}(\textbf{u}):= \sum_{i=1}^{N}\mathcal{J}_{i}(\textbf{u}), \quad \U^N:=\prod_{i=1}^{N}\U_{i}.
\end{aligned}\right.
\end{equation}MG, MT are both \emph{decentralized} on $L_{\mathbb{F}_i}^2(0,T;\R^m)$ or $\prod_{i=1}^{N} L_{\mathbb{F}_i}^2(0,T;\R^m)$. So, they differ from mean-field-type control (\cite{ad2010}, \cite{p2016}, etc.), still centralized on $\mathbb{F}.$ But, MT aims to unified \emph{team cost} $\mathcal{J}_{soc}^{(N)}$ by \emph{cooperative} team-decision $\textbf{u}=(u_1, \cdots, u_{N})$, whereas MG to \emph{competitive} $\mathcal{J}_{i}$ by individual $u_{i}$. Our inspiration is to recast MG, MT from a \emph{coalition matrix} insight.

\subsection{A coalition matrix representation} \label{coa}

First, we introduce \emph{nominal cost vector}
$\mathbf{J} := (\mathcal{J}_1, \cdots, \mathcal{J}_{N})^{T}$ formalized by all $\{\mathcal{J}_i\}_{i=1}^{N}.$ For fixed $\mathcal{A}_{i},$ we call $\mathcal{J}_{i}$  its \emph{principal} cost, while $\{\mathcal{J}_{k}\}_{k \neq i}$ \emph{marginal} costs. Next, set \emph{effective cost vector} $\mathbf{K}:=(\mathcal{K}_1, \cdots, \mathcal{K}_{N})^{T}$ with $\mathcal{K}_{i}$ the effective cost really ``targeted" by $\mathcal{A}_{i}.$
In general, $\mathcal{K}_{i} \neq \mathcal{J}_{i}.$

By matrix operation, $\mathbf{K}=\mathbf{C} \cdot \mathbf{J}$ for some \emph{coalition matrix} $\mathbf{C}=[c_{i,j}]_{1 \leq i, j \leq N}$ for which some typical forms are listed below as $\mathbf{C}_1$-$\mathbf{C}_5$:

\begin{figure}[htbp]
\centering
\vspace*{0cm}
\includegraphics[height=2cm,width=8.5cm]{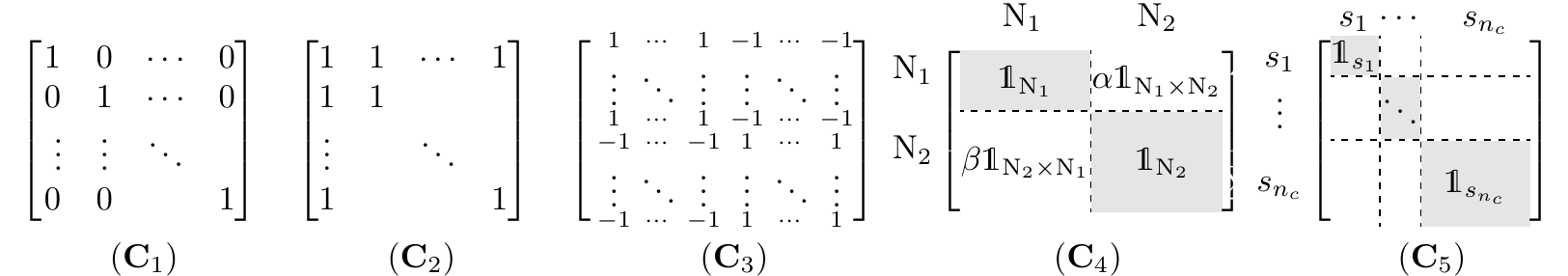}
\vspace*{0cm}
\caption{$\mathbf{C}-$form coalition matrices.}\label{pic1}
\end{figure}

Notably, all coalition $\mathbf{C}-$forms above share a common \emph{block} feature: the diagonal
``covered" by $1\leq n_c\leq N$ square \textbf{one} matrix $\mathds{1}_{s}$ (i.e., all entries are $1$) of possible
varying sizes (e.g., $s_1,\cdots,s_{n_c}$) as exhibited by the general $\mathbf{C}_5$. Noting $\mathbf{C}_5$ nests $\mathbf{C}_4$ if
$n_c=2,\ s_1 = N_1,\  s_2 = N_2 = N-N_1$ and properly setting its off-diagonal blocks; $\mathbf{C}_4$ further nests $\mathbf{C}_3$ by $\alpha=\beta=-1$ and $\textbf{C}_2$ by $\alpha=\beta=1$ $(n_c = 1,\ s_1 = N)$; also, $\mathbf{C}_5$ nests $\mathbf{C}_1$ by $n_c=N,\ s_i\equiv 1$. The meaning of ``coalition" is illustrated by classifications below:
\begin{equation}\label{cform}\left\{\begin{aligned}
\mathbf{C}&=\mathbf{C}_{1}=\textbf{I}_{N}\Longrightarrow (\mathbf{K}=\mathbf{J})\Leftrightarrow (\mathcal{K}_{i}=\mathcal{J}_{i})\\
 &\Longrightarrow \text{\textbf{MG}}\Leftrightarrow (\text{non-cooperative})\ \{\mathcal{A}_{i}\}_{i=1}^{N};\ \\
 \mathbf{C}&=\mathbf{C}_{2}=\mathds{1}_{N} \Longrightarrow (\mathcal{K}_{i}\equiv\mathcal{J}^{(N)}_{soc}=\sum \mathcal{J}_{k})\\
  &\Longrightarrow \text{\textbf{MT}}\Leftrightarrow (\text{cooperative})\ \{\mathcal{A}_{i}\}_{i=1}^{N}.\end{aligned}\right.\end{equation}
So, different $\mathbf{C}-$forms induce different \emph{coalitions} among $\{\mathcal{A}_{i}\}_{i=1}^{N}.$ Here, we do not distinguish more subtle \emph{exogenous} or \emph{endogenous} coalition formation (e.g., \cite{MTMW09}, \cite{STY08}, \cite{Z94}) as they are economics biased. We observe $\mathbf{C}_{1}, \mathbf{C}_{2}$ as exactly two extreme $\mathbf{C}-$forms which are leading to ``pure" MG or MT. But various other $\mathbf{C}-$forms $1<n_c<N$ do exist such as $\mathbf{C}_{3}\-\mathbf{C}_{5}$ which should connect to other meaningful structures beyond ``pure" game/team. In fact, by setting $\I_1=\{\theta_1,\cdots,\theta_{N_1}\}$, $\I_2=\{\vartheta_1,\cdots,\vartheta_{N_2}\}$ and $N_1+N_2=N$, where $\theta_i,\ \vartheta_j$ are all positive indexes for $1\leq i\leq N_1,\ 1\leq j\leq N_2$, we may proceed more along typical $\textbf{C}_4:$
\begin{equation}\label{c4}\begin{aligned}
&\mathbf{C}=\mathbf{C}_{4}  \Longrightarrow \mathcal{K}_{i}=\\
&\left\{\begin{aligned}& \mathbf{J}_{\text{mix}}^{(1)}:= \mathcal{J}^{1,(N)}_{\text{soc}} \+ \alpha \mathcal{J}^{2,(N)}_{\text{soc}},\quad \mathcal{J}^{1,(N)}_{\text{soc}}:= \sum_{k \in \I_{1}} \mathcal{J}_{k}, \\ &
\mathbf{J}_{\text{mix}}^{(2)}:= \mathcal{J}^{2,(N)}_{\text{soc}} \+ \beta \mathcal{J}^{1,(N)}_{\text{soc}},\quad  \mathcal{J}^{2,(N)}_{\text{soc}}:=\sum_{k \in \I_{2}} \mathcal{J}_{k},
\end{aligned}\right.\end{aligned}\end{equation}
where $\alpha,\ \beta\in \R$. For $k=1,2$, let $\pi_k^{(N)}=\frac{N_k}{N}$, then $\pi^{(N)}=(\pi_1^{(N)},\pi_2^{(N)})$ is a probability vector representing the empirical distribution of $\I_1$ and $\I_2$.
Thus, the  original LS system divides into two sub-systems: LS1 $:= \{\mathcal{A}_{k}\}_{k\in \I_{1}}$, LS2 $:= \{\mathcal{A}_{k}\}_{k\in \I_{2}}$. All agents inside LS1/LS2 aim the same cost so formalize two cooperative \emph{internal} MTs. LS1 and LS2 remain competitive outside if $\alpha<0,\ \beta<0$; while cooperative outside if $\alpha>0,\ \beta>0$.

\subsection{Mixed game-team via coalition}
Now, we are ready to formally introduce the ``mixed" game-team (\textbf{Mix}) via coalition matrix representation $\mathbf{C}_{4}$ as
\begin{equation}\label{mix}
\text{\textbf{(Mix)}}
\left\{\begin{aligned}
&\text{\textbf{LS1:}} \inf_{\textbf{u}_{1}(\cdot) \in \mathbf{U}^{({1})}}\mathbf{J}_{\text{mix}}^{(1)}(\textbf{u}_{1}; \textbf{u}_{2}), \\ &\quad\textbf{u}_{1}:=(\cdots, u_{k}, \cdots)|_{k \in \I_{1}},\ \mathbf{U}^{(1)}:=\prod_{k \in \I_{1}}\U_{k};\\
 &\text{\textbf{LS2:}} \inf_{\textbf{u}_{2}(\cdot) \in \mathbf{U}^{({2})}}\mathbf{J}_{\text{mix}}^{(2)}(\textbf{u}_{2}; \textbf{u}_{1}),\\
 &\quad \textbf{u}_{2}:=(\cdots, u_{k}, \cdots)|_{k \in \I_{2}}, \ \mathbf{U}^{(2)}:=\prod_{k \in \I_{2}}\U_{k}.
 \end{aligned}\right.
 \end{equation}
To better illustrate the categories of \textbf{LS1} and \textbf{LS2}, we need rewrite the dynamics $\mathcal{S}_i$ in \eqref{ls} as the following homogeneous form:
\begin{equation}\label{state} \left\{\begin{aligned}
    &dx_i(t)\= \left(A_1x_i(t)\+B_1u_i(t)\+ F_1x^{(N)}(t) \right)dt\\
    &\qquad\qquad\+\sigma_1(t)dW_{i}(t), \quad x_i(0)\=\xi_{i},\quad i\in \I_1,\\
    &dx_j(t)\= \left(A_2x_j(t)\+B_2u_j(t)\+ F_2x^{(N)}(t)\right )dt\\
    &\qquad\qquad\+\sigma_2(t)dW_{j}(t), \quad x_j(0)\=\eta_{j},\quad j\in \I_2,
\end{aligned}\right.\end{equation}
where $\xN (\cdot)= \frac{1}{N}\left(\sum_{i\in\I_1}x_i(\cdot) \+ \sum_{j\in\I_2}x_j(\cdot)\right)$. For the sake of notation simplicity, we denote  $\bx_1 := (x_{\theta_1},\cdots,x_{\theta_{N_1}})$, $\bx_2 := (x_{\vartheta_1},\cdots,x_{\vartheta_{N_2}})$. we also rewrite the individual cost functionals \eqref{JMG} as
\s\begin{equation} \label{cost}
\left\{\begin{aligned}
&\mathcal{J}_{i}(\bu_1(\cdot),\bu_2(\cdot))\\
&=\frac{1}{2}\mathbb{E}\int_0^{T}\Big[\|x_{i}(t)-\Gamma_1 x^{(N)}(t)\|_{Q_1}^{2}\+\|u_{i}(t)\|^{2}_{R_1}\Big]dt,\\
&\mathcal{J}_{j}(\bu_1(\cdot),\bu_2(\cdot))\\
&=\frac{1}{2}\mathbb{E}\int_0^{T}\Big[\|x_{j}(t)-\Gamma_2 x^{(N)}(t)\|_{Q_2}^{2}\+\|u_{j}(t)\|^{2}_{R_2}\Big]dt.
\end{aligned}\right.
\end{equation}\n

The agents in $\{\A_{i}\}_{i\in\I_1}$ and $\{\A_{j}\}_{j\in\I_2}$ are cooperative respectively and their social cost functionals $\mathbf{J}_{\text{mix}}^{(1)}$, $\mathbf{J}_{\text{mix}}^{(2)}$ are given by \eqref{c4}. We impose the following general assumptions, {which  are commonly used in LQ models}, on the coefficients:
\begin{description}
  \item[(H1)] $A_k,\ F_k,\ \Gamma_k\in \R^{n\times n}$, $B_k\in \R^{n\times m}$, $\sigma_k\in L^\infty(0,T;\mathbb R^{n})$, $k=1,2$.
  \item[(H2)] $Q_k\in \S^n$, $R_k\in \S^m$, $Q_k> 0,\ R_k> 0,$ $k=1,2$.
  \item[(H3)] $\{\xi_i\}_{i=1}^{N_1}$ are independent identically distributed (i.i.d) with mathematical expectation $\E\xi$; $\{\eta_j\}_{j=1}^{N_2}$ are i.i.d with mathematical expectation $\E\eta$.
  \item[(H4)] There exists a probability mass vector $\pi=(\pi_1, \pi_2)$ such that $\displaystyle{\lim_{N\rightarrow+\infty}}\pi^{(N)}=\pi$, $\displaystyle{\min_{1 \leq k \leq 2}}\pi_{k}>0.$
\end{description}
Thus, we can propose the following homogeneous mixed game-team  problem:
\begin{Problem}\label{problem1}
  Find a centralized strategy set $\left(\bar{\bu}_1,\bar{\bu}_2\right)$, where $\bar{\bu}_1\! =\! (\bar{u}_{\theta_1},\!\cdots\!,\bar{u}_{\theta_{N_1}})$, $\bar{\bu}_2 = (\bar{u}_{\vartheta_1},\cdots,\bar{u}_{\vartheta_{N_2}})$, $\bar u_i\in \mathcal U_i$, $\bar u_j\in \mathcal U_j$, $i\in \I_1$, $j\in \I_2$, such that
  \begin{equation}\label{mix2}
  \emph{\text{\textbf{(Mix)}}}
  \left\{\begin{aligned}
  &\text{\rm\textbf{LS1:}}\quad \mathbf{J}_{\text{\rm mix}}^{(1)}(\bar{\bu}_1,\bar{\bu}_2) = \inf_{{\bu}_1\in \mathbf{U}^{({1})}} \mathbf{J}_{\text{\rm mix}}^{(1)}( {\bu}_1,\bar{\bu}_2),\\
  &\text{\rm\textbf{LS2:}}\quad \mathbf{J}_{\text{\rm mix}}^{(2)}(\bar{\bu}_1,\bar{\bu}_2) = \inf_{{\bu}_2 \in \mathbf{U}^{({2})}} \mathbf{J}_{\text{\rm mix}}^{(2)}( \bar{\bu}_1, {\bu}_2).
  \end{aligned}\right.
  \end{equation}
\end{Problem}

\textbf{(Mix)} is meaningful, not only in mathematics as verified by coalition arguments above, but also in practice as supported by real motivations (see \cite{BC98}, \cite{DHM20}, etc).

\section{Variation synthesis analysis}\label{sec 3}
Initially, we synthesize all response components for $\mathbf{J}_{\text{mix}}^{(1)}$ in \eqref{mix2}, from a bilateral viewpoint of \textbf{LS1}, \textbf{LS2}. In what follows, we only focus on the viewpoint of \textbf{LS1}, and the similar argument can be applied to \textbf{LS2}.

\subsection{Variation decomposition}
Let $\bar{\bu}_1$, $\bar{\bu}_2$ be
centralized optimal strategies of the agents in $\I_1$, $\I_2$. We now perturb
$u_i$ and keep $\bar u_{-i}=(\bar u_{\theta_1}, \cdots, \bar u_{i-1}, \bar u_{i+1},
\cdots, \bar u_{\theta_{N_1}})$, $\bar{\bu}_2$ fixed. For $k\neq i$, denote the
perturbation $\d u_i=u_i-\bar u_i$, $\d x_i=x_i-\bar x_i$, $\d
x_k=x_k-\bar x_k$, $\d
x^{(N)}=\frac{1}{N}\sum_{k\in\I_1\cup\I_2} \d x_k$,  and
$\d\mathcal{J}^{1,(N)}_{\text{soc}}$, $\d\mathcal{J}^{2,(N)}_{\text{soc}}$ are the first variations (Fr\'{e}chet differentials)
of $\mathcal{J}^{1,(N)}_{\text{soc}}$, $\mathcal{J}^{2,(N)}_{\text{soc}}$ w.r.t.  $\d u_i$. Therefore, 
\begin{equation*}
\left\{\begin{aligned}
&d\d x_i(\d u_i)\= \left[A_1\d x_i(\d u_i)\+B_1\d u_i\+ F_1\d x^{(N)} (\d u_i) \right]dt, \\
&\quad\  \d x_i(0)\=0,\\
&d\d x_k(\d u_i)\= \left[A_1\d x_k(\d u_i) \+ F_1\d x^{(N)}(\d u_i) \right]dt, \\
&\quad\ \d x_k(0)\= 0,\  k\in \I_1,\  k\neq i,\\
&d\d x_k(\d u_i)\= \left[A_2\d x_k(\d u_i) \+ F_2\d x^{(N)}(\d u_i) \right]dt, \\
&\quad\ \d x_k(0)\= 0,\  k\in \I_2.\\
\end{aligned}\right.
\end{equation*}
Then we derive the following lemma.
\begin{Lemma}
$\d \mathbf{J}_{\text{mix}}^{(1)}(\d u_{i})$ can be represented as
\s\begin{equation}\label{3.5}
 \begin{aligned}
&\d \mathbf{J}_{\text{mix}}^{(1)}(\d u_{i}) =\mathbb{E}\int_0^{T}\Big[ \left\langle\bar{\Theta}_{1}, \d x_{i}\right\rangle
\+\left\langle\Theta_{2}, \d u_{i}\right\rangle\+\Big\langle\bar{\Theta}_{3}, \sum_{ k\in\I_1,k\neq i}\d x_{k} \Big\rangle
\\
&\ +\Big\langle\bar{\Theta}_{4}, \sum_{k \in \I_{2}}\d x_{k}\Big\rangle\+\sum_{k\in\I_1,k\neq i}\left\langle\Theta_{5}^{k}, \d x_{k}\right\rangle
\+\sum_{k \in \I_{2}}\left\langle\Theta_{6}^{k}, \d x_{k}\right\rangle\Big] dt \+\varepsilon^i_1, \\
\end{aligned}
\end{equation}\n
where
\s\begin{equation*}
\begin{aligned}
\varepsilon^i_1 = \mathbb{E}\int_0^{T}\Big[ &\left\langle\Theta_1 - \bar{\Theta}_1,\d x_{i}\right\rangle \+ \Big\langle\Theta_3 - \bar{\Theta}_3,\sum_{k\in\I_1,k\neq i}\d x_{k}\Big\rangle\\
 &\+ \Big\langle{\Theta}_{4} - \bar{\Theta}_{4}, \sum_{k \in \I_{2}}\d x_{k} \Big\rangle\Big] dt,
\end{aligned}
\end{equation*}\n
\tiny\begin{equation*}
\left\{\begin{aligned}
     & \Theta_1 =  Q_1\x_i \- \left( Q_1  \Gamma_1 \x^{(N)} \+ \Gamma_1^TQ_1  \pi_1^{(N)}\frac{\sum_{k\in \I_1}\x_{k}}{N_1} \- \pi_1^{(N)} \Gamma_1^TQ_1 \Gamma_1 \x^{(N)}\right) \\
     &\qquad\qquad\- \alpha\left( \Gamma_2^TQ_2  \pi_2^{(N)}\frac{\sum_{k\in \I_2}\x_{k} }{N_2} \- \pi_2^{(N)} \Gamma_2^TQ_2 \Gamma_2 \x^{(N)}\right),\\
     & \Theta_3 =- \left( Q_1  \Gamma_1 \x^{(N)} \+ \Gamma_1^TQ_1  \pi_1^{(N)}\frac{\sum_{k\in \I_1}\x_{k}}{N_1} \- \pi_1^{(N)} \Gamma_1^TQ_1 \Gamma_1 \x^{(N)}\right) \\
          &\qquad\qquad\- \alpha\left( \Gamma_2^TQ_2  \pi_2^{(N)}\frac{\sum_{k\in \I_2}\x_{k} }{N_2} \- \pi_2^{(N)} \Gamma_2^TQ_2 \Gamma_2 \x^{(N)}\right),   \\
 & \Theta_4 = - \alpha\left(Q_2  \Gamma_2 \x^{(N)} \+ \Gamma_2^TQ_2  \pi_2^{(N)}\frac{\sum_{k\in \I_2}\x_{k}}{N_2} \- \pi_2^{(N)} \Gamma_2^TQ_2 \Gamma_2 \x^{(N)} \right) \\
     &\qquad\qquad\- \left( \Gamma_1^TQ_1  \pi_1^{(N)}\frac{\sum_{k\in \I_1}\x_{k} }{N_1}\- \pi_1^{(N)} \Gamma_1^TQ_1 \Gamma_1 \x^{(N)} \right),            \\
     & \Theta_2 = R_1, \quad  \Theta_5^{k} = Q_1\bar{x}_{k}, \quad  \Theta_6^{k} = \alpha Q_2\bar{x}_{k}, \\
 &\bar{\Theta}_1 = Q_1\x_i \- \left( Q_1  \Gamma_1 \m \+ \Gamma_1^TQ_1  \pi_1\m_1 \- \!\pi_1 \Gamma_1^TQ_1 \Gamma_1 \m\right)\\
       &\qquad\qquad \- \alpha\left( \Gamma_2^TQ_2  \pi_2\m_2 \- \!\pi_2 \Gamma_2^TQ_2 \Gamma_2 \m\right),\\
      &\bar{\Theta}_3 =- \left( Q_1  \Gamma_1 \m \+ \Gamma_1^TQ_1  \pi_1\m_1 \- \pi_1 \Gamma_1^TQ_1 \Gamma_1 \m\right) \\
      &\qquad\qquad \- \alpha\left( \Gamma_2^TQ_2  \pi_2\m_2 \- \pi_2 \Gamma_2^TQ_2 \Gamma_2 \m\right),    \\
     \end{aligned}\right.
\end{equation*}
\begin{equation*}
\begin{aligned}
         &\bar{\Theta}_4 = - \alpha\left(Q_2  \Gamma_2 \m \+ \Gamma_2^TQ_2  \pi_2\m_2 \- \pi_2 \Gamma_2^TQ_2 \Gamma_2 \m \right)\\
     &\qquad\qquad \- \left( \Gamma_1^TQ_1  \pi_1\m_1\- \pi_1 \Gamma_1^TQ_1 \Gamma_1 \m \right).           \\
 \end{aligned}
\end{equation*}
\n
Here, $\m_1$, $\m_2$ and $\m = \pi_1\m_1 \+ \pi_2\m_2$ are mean-field (MF) approximations of $\frac{1}{N_1}\sum_{k\in \I_1}\x_{k}$, $\frac{1}{N_2}\sum_{k\in \I_2}\x_{k}$ and  $\bar{x}^{(N)}$, respectively.
\end{Lemma}

\subsection{Bilateral duality}
\begin{Lemma}
$\d \mathbf{J}_{\text{mix}}^{(1)}(\d u_{i})$ can further be represented as
\s\begin{equation}\label{2.3}
 \begin{aligned}
&\d \mathbf{J}_{\text{mix}}^{(1)}(\d u_{i}) = \mathbb{E}\int_0^{T} \Big[\left\langle\bar{\Theta}_{1} \+ \pi_1  F_1^T\E p_k^{(1)*}  \+ \pi_2  F_2^T\E p_k^{(2)*}\right.\\
&\qquad\left.\+  \pi_1F_1^Tp^{(1)*} \+ \pi_2F_2^Tp^{(2)*}, \d x_{i}\right\rangle+\left\langle\Theta_{2}  , \d u_{i}\right\rangle\ \Big] dt\\
&\qquad \+\varepsilon^i_1  \+ \varepsilon^i_2 \+ \varepsilon^i_3, \\
\end{aligned}
\end{equation}\n
where
\scriptsize\begin{equation*}
\begin{aligned}
&\varepsilon^i_2 \= \mathbb{E}\int_0^{T} \Big[\Big\langle\bar{\Theta}_{3},  \sum_{k\in\I_1 ,k\neq i}\d x_{k} \-  x^{(1),**}\Big\rangle+\Big\langle\bar{\Theta}_{4},  \sum_{k\in\I_2}\d x_{k} \-  x^{(2),**}\Big\rangle \\
&\qquad \+\frac{1}{N_1}\sum_{k \in \I_{1},k\neq i}\left\langle\Theta_{5}^{k}, N_1\d x_k \-  x^{(1),*}_k\right\rangle\\
&\qquad+\frac{1}{N_2}\sum_{k \in \I_{2}}\left\langle\Theta_{6}^{k},   N_2\d x_k \- x^{(2),*}_k\right\rangle \Big]  dt,\\
&\varepsilon^i_3 \= \mathbb{E}\int_0^{T} \Big\langle  \pi_1  F_1^T \Big(\frac{  \sumI1 p_k^{(1)}  }{N_1} \- \E p_k^{(1)*}\Big)\\
 &\qquad\+ \pi_2F_2^T\Big( \frac{\sum_{k\in\I_2 }p_k^{(2)}   }{N_2} \- \E p_k^{(2)*} \Big)\+\pi_1F_1^T\left(p^{(1)} \- p^{(1)*}\right)  \\
&\qquad \+ \pi_2F_2^T\left(p^{(2)} \- p^{(2)*}\right), \d x_{i}\Big\rangle dt,
\end{aligned}
\end{equation*}\n
and
\scriptsize\begin{equation*}
  \left\{\begin{aligned}
      &d   x^{(1),*}_{k} \=\Big[A_1  x^{(1),*}_{k} \+ \pi_{1} F_1(\d x_i \+ x^{(1),**} \+  x^{(2),**})\Big]dt, \\
      &d  x^{(2),*}_{k} \=\Big[A_2  x^{(2),*}_{k} \+ \pi_{2} F_2(\d x_i \+  x^{(1),**} \+  x^{(2),**})\Big]dt,\\
     & d x^{(1),**} \= \Big[A_1   x^{(1),**} \+ \pi_{1}F_1(  x^{(1),**} \+  x^{(2),**} \+ \d x_i)\Big]dt,\\
     &d x^{(2),**} \= \Big[A_2  x^{(2),**} \+ \pi_{2}F_2(  x^{(1),**} \+  x^{(2),**} \+ \d x_i)\Big]dt,  \\
     & x^{(1),*}_{k}(0)=0, \  \ x^{(2),*}_{k}(0)=0, \ x^{(1),**}(0)=0, \ x^{(2),**}(0) \= 0,
  \end{aligned}\right.
\end{equation*}\n
\tiny\begin{equation*}
  \left\{\begin{aligned}
     &  dp_k^{(1)} = - \left(\Theta_{5}^{k} \+ A_1^T p_k^{(1)} \right) dt \+ q_k^{(1)} dW_k \+  {\sum_{k'\neq k}} q_{kk'}^{(1)} dW_{k'},\\
     & dp_k^{(2)} = -\left( \Theta_{6}^{k} \+ A_2^T p_k^{(2)} \right) dt \+ {q}_k^{(2)} dW_k \+  {\sum_{k'\neq k} }q_{kk'}^{(2)} dW_{k'},\\
     &   dp^{(1)} = - \Big(\bar{\Theta}_{3} \+ \pi_1  {\frac{\sumI1 F_1^Tp_k^{(1)} }{N_1}} \+ \pi_2  {\frac{\sum_{k\in\I_2 } F_2^Tp_k^{(2)} }{N_2}} \\
     &\hspace{1.8cm}\+ A_1^Tp^{(1)}  \+  \pi_1F_1^Tp^{(1)} \+ \pi_2F_2^Tp^{(2)}  \Big) dt\+  {\sum} q_{k'}^{(1)} dW_{k'}, \\
     &  d{p}^{(2)} = -\Big( \bar{\Theta}_{4} \+ \pi_1  {\frac{\sumI1 F_1^Tp_k^{(1)} }{N_1}} \+ \pi_2  {\frac{\sum_{k\in\I_2 } F_2^Tp_k^{(2)} }{N_2}}\\
     &\hspace{1.8cm} \+  \pi_1F_1^Tp^{(1)} \+ \pi_2F_2^Tp^{(2)}\+ A_2^Tp^{(2)}  \Big) dt \+  {\sum} q_{k'}^{(2)} dW_{k'}, \\
& p_k^{(1)}(T)=0,\ p_k^{(2)}(T)=0,\ p^{(1)}(T)=0,\ {p}^{(2)}(T) = 0,
  \end{aligned}\right.
\end{equation*}\n
\scriptsize\begin{equation*}
  \left\{\begin{aligned}
     &\widehat{\mathcal{P}}^{(2)*}_k: dp_k^{(1)*}=  - \left(\Theta_{5}^{k} \+ A_1^T p_k^{(1)*}\right) dt \+ q_k^{(1)*} dW_k
      ,\\
     &\widehat{\mathcal{P}}^{(1)*}_k: dp_k^{(2)*} = -\left( \Theta_{6}^{k} \+ A_2^T p_k^{(2)*}   \right) dt \+ q_k^{(2)*} dW_k,\\
     &\widehat{\mathcal{P}}^{(2)*}:  dp^{(1)*}   \=  - \left[\bar{\Theta}_{3} \+ \pi_1   F_1^T\E p_k^{(1)*}   \+ \pi_2   F_2^T\E p_k^{(2)*} \right.\\
      &\qquad\qquad\left.\+ \left(A_1^T   \+  \pi_1F_1^T \right)p^{(1)*} \+ \pi_2F_2^Tp^{(2)*}  \right] dt, \\
     &\widehat{\mathcal{P}}^{(1)*}: d{p}^{(2)*} = -\left[ \bar{\Theta}_{4} \+ \pi_1    F_1^T\E p_k^{(1)*}   \+ \pi_2   F_2^T\E p_k^{(2)*} \right.\\
      &\qquad\qquad\left.\+  \pi_1F_1^Tp^{(1)*} \+ \left(\pi_2F_2^T \+ A_2^T \right)p^{(2)*} \right] dt,\\
     &  p_k^{(1)*}(T) =0,\    p_k^{(2)*}(T) =0,\   p^{(1)*}(T)  =0,\    {p}^{(2)*}(T)= 0.
  \end{aligned}\right.\\
\end{equation*}\n
\end{Lemma}

An asymptotic ``Fr\'{e}chet response" holds for $\delta J_{i}^{(1)}=\lim  \d \mathbf{J}_{\text{mix}}^{(1)}(\delta u_{i})$:
 \begin{equation}\label{fre}
 \begin{aligned}&\text{\textbf{LS1}:} \ \ \delta J_{i}^{(1)} \!=\!\langle{\Theta}^{\dag}(\bar{x}_{i}, \textbf{m}, \textbf{P}^*), \delta x_{i}\rangle\!+\!\langle\Theta^{\dag\dag}(\bar{u}_{i}), \delta u_{i}\rangle \end{aligned}
 \end{equation}
 with $\Theta^{\dag\dag}=R_1 \bar{u}_{i}$ and $\Theta^{\dag}$ is revised on $\left(\bar{\Theta}_1,\bar{\Theta}_3,\bar{\Theta}_4, {\Theta}_5^{k}, {\Theta}_6^{k}\right)$. Here, $\textbf{m}:=(\m, \m_1, \m_2)$, $\widehat{\textbf{P}}^*:=(\widehat{\mathcal{P}}_{k}^{(1)*}, \widehat{\mathcal{P}}^{(1)*}; \widehat{\mathcal{P}}_{k}^{(2)*}, \widehat{\mathcal{P}}^{(2)*})$.

\section{Distributed design} \label{sec 4}

\subsection{Auxiliary control} \label{auxiliary}
Motivated by \eqref{2.3}, we introduce the following auxiliary problem for $i\in\I_1$:
\begin{Problem} \label{problem2}
Minimize $J_i(u_i)$ over $u_i\in\mathcal U_i$ where
\s\begin{equation*}
\left\{\begin{aligned}
     &dx_i\= (A_1x_i\+B_1u_i\+ F_1 \m )dt\+\sigma_1dW_{i}, \  x_i(0)\=\xi_{i},\  i\in \I_1,                            \\
     & J_i^{(1)}(u_i) \= \frac{1}{2}\E  \int_0^T [\langle Q_1 x_i, x_i\rangle\+ 2\langle S_1, x_i\rangle\+\langle R_1 u_i, u_i\rangle] dt,\\
     &S_1=  \- \left( Q_1  \Gamma_1 \m \+ \Gamma_1^TQ_1  \pi_1\m_1 \- \pi_1 \Gamma_1^TQ_1 \Gamma_1 \m\right)\\
     &\hspace{1cm} \- \alpha\left( \Gamma_2^TQ_2  \pi_2\m_2 \- \pi_2 \Gamma_2^TQ_2 \Gamma_2 \m\right)+ \pi_1  F_1^T\E p_k^{(1)*} \\
     &\hspace{1cm} \+ \pi_2  F_2^T\E p_k^{(2)*} \+  \pi_1F_1^Tp^{(1)*} \+ \pi_2F_2^Tp^{(2)*}.
  \end{aligned}\right.
  \end{equation*}\n
\end{Problem}
The mean-field terms $\m$, $\m_1$, $\m_2$, $p_k^{(1)*}$, $ p_k^{(2)*}$ will be determined by the CC system later. Similarly, for $j\in\I_2$, by following the procedure in Section \ref{sec 3},  we can also  introduce the following auxiliary problem
\begin{Problem} \label{problem 3}
Minimize $J_j(u_j)$ over $u_j\in\mathcal U_j$ where
\s \begin{equation*}
\left\{\begin{aligned}
     &dx_j\= (A_2x_j\+B_2u_j\+ F_2 \m )dt\+\sigma_2dW_{j}, \  x_j(0)\=\eta_{j},\  j\in \I_2,                            \\
     & J_j^{(2)}(u_j) \= \frac{1}{2}\E  \int_0^T [\langle Q_2 x_j, x_j\rangle\+ 2\langle S_2, x_j\rangle\+\langle R_2 u_j, u_j\rangle] dt,\\
     &S_2=  \- \left( Q_2  \Gamma_2 \m \+ \Gamma_2^TQ_2  \pi_2\m_2 \- \pi_2 \Gamma_2^TQ_2 \Gamma_2 \m\right)\\
     &\hspace{1cm} \- \alpha\left( \Gamma_1^TQ_1  \pi_1\m_1 \- \pi_1 \Gamma_1^TQ_1 \Gamma_1 \m\right)+ \pi_2  F_2^T\E \hat{p}_k^{(2)*}\\
     &\hspace{1cm} \+ \pi_1  F_1^T\E \hat{p}_k^{(1)*}  \+  \pi_2F_2^T\hat{p}^{(2)*} \+ \pi_1F_1^T\hat{p}^{(1)*},
  \end{aligned}\right.
  \end{equation*}\n
\end{Problem}
where the limiting duality on $\delta u_{j}$ yields parallel $\widehat{\textbf{P}}^*:=(\widehat{\mathcal{P}}_{k}^{(1)*}, \widehat{\mathcal{P}}^{(1)*}; \widehat{\mathcal{P}}_{k}^{(2)*}, \widehat{\mathcal{P}}^{(2)*})$.
Here, it holds that
\scriptsize\begin{equation*}
  \left\{\begin{aligned}
  &\widehat{\mathcal{P}}^{(2)*}_k:   d\hat{p}^{(2)*}_k =  - \left(\hat{\Theta}_{5}^{k} \+ A_2^T \hat{p}_k^{(2)*}\right) dt \+ \hat{q}_k^{(2)*} dW_k
      ,                                                                     \\
     & \widehat{\mathcal{P}}^{(1)*}_k:  d\hat{p}^{(1)*}_k = -\left( \hat{\Theta}_{6}^{k} \+ A_1^T \hat{p}_k^{(1)*}  \right) dt \+ \hat{q}_k^{(1)*} dW_k
     ,                                                                  \\
     &\widehat{\mathcal{P}}^{(2)*}:   d\hat{p}^{(2)*} =  - \left[\hat{\bar{\Theta}}_{4} \+ \pi_2  F_2^T\E \hat{p}_k^{(2)*}   \+ \pi_1   F_1^T\E \hat{p}_k^{(1)*}\right.\\
     &\left.\qquad\qquad\+ \left(A_2^T   \+  \pi_2F_2^T \right)\hat{p}^{(2)*} \+ \pi_1F_1^T\hat{p}^{(1)*}  \right] dt,  \\
     & \widehat{\mathcal{P}}^{(1)*}: d\hat{p}^{(1)*} = -\left[ \hat{\bar{\Theta}}_{3} \+ \pi_2   F_2^T\E \hat{p}_k^{(2)*}   \+ \pi_1   F_1^T\E \hat{p}_k^{(1)*} \right.\\
       &\left.\qquad\qquad\+  \pi_2F_2^T \hat{p}^{(2)*} \+ \left(\pi_1F_1^T \+ A_1^T \right)\hat{p}^{(1)*} \right] dt,\\
          & \hat{p}_k^{(2)*}(T) = \hat{p}_k^{(1)*}(T) = \hat{p}^{(2)*}(T) = \hat{p}^{(1)*}(T) = 0,
  \end{aligned}\right.\\
\end{equation*}\n
and
\scriptsize\begin{equation*}
  \left\{\begin{aligned}
      &\hat{\bar{\Theta}}_4 =- \left( Q_2  \Gamma_2 \m \+ \Gamma_2^TQ_2  \pi_2\m_2 \- \pi_2 \Gamma_2^TQ_2 \Gamma_2 \m\right)\\
      &\qquad\qquad\- \beta\left( \Gamma_1^TQ_1  \pi_1\m_1 \- \pi_1 \Gamma_1^TQ_1 \Gamma_1 \m\right),                   \\
     & \hat{\bar{\Theta}}_3 =- \beta\left(Q_1  \Gamma_1 \m \+ \Gamma_1^TQ_1  \pi_1\m_1 \- \pi_1 \Gamma_1^TQ_1 \Gamma_1 \m \right) \\
      &\qquad\qquad\- \left( \Gamma_2^TQ_2  \pi_2\m_2\- \pi_2 \Gamma_2^TQ_2 \Gamma_2 \m \right),        \\
     &\hat{\Theta}_{5}^{k} = Q_2\bar{x}_{k},\quad \hat{\Theta}_{6}^{k} = \beta Q_1\bar{x}_{k}.
  \end{aligned}\right.
\end{equation*}\n
 The above analysis constructs a bilateral auxiliary control problem
\begin{equation}\label{ac}\begin{aligned}
&\text{\textbf{(Bilateral auxiliary problem)}}\\
&
\left\{
\begin{aligned}& \mathcal{A}_{i}:\inf_{u_{i}(\cdot) \in \U_{i}}J_{i}^{(1)}(u_{i}; \textbf{m}, \textbf{P}^*), \\& {\mathcal{A}}_{j}: \inf_{u_{j}(\cdot) \in \U_{j}}J_{j}^{(2)}(u_{j}; \textbf{m}, \widehat{\textbf{P}}^*).
\end{aligned}
\right.
\end{aligned}\end{equation}
Then, upon distributed $\U_{i}$ or $\U_{j}$, two generic $\mathcal{A}_{i}, \mathcal{A}_{j}$ aim to optimize $J_{i}^{(1)}$ or $J_{j}^{(2)}$. Applying standard LQ method, under some mild (say, convexity) condition, above bilateral auxiliary problem can be solved with optimal pair $(\check{x}_{i}(\check{u}_{i}), \check{u}_{i})$ for $\mathcal{A}_{i},$ $(\check{x}_{j}(\check{u}_{j}), \check{u}_{j})$ for $\mathcal{A}_{j},$ which depend on $(\textbf{m}, \textbf{P}^*)$ or $(\textbf{m}, \widehat{\textbf{P}}^*)$ by \eqref{ac}. Note that $\textbf{m}$ is not specified yet.
By stochastic maximum principle, we have the following result for the bilateral auxiliary problem:
\begin{Proposition}
  Under (H1)-(H4), Problem \ref{problem2}-\ref{problem 3} are both uniquely solvable, and the optimal auxiliary controls $\u_i$, $\u_j$ are determined by the following Hamiltonian systems
   \begin{equation*}
   \left\{\begin{aligned}
   &d\cx_i\= (A_1\cx_i\+B_1\cu_i\+ F_1 \m )dt\+\sigma_1dW_{i}, \quad \cx_i(0)\=\xi_{i}, \\
   &d\cx_j\= (A_2\cx_j\+B_2\cu_j\+ F_2 \m )dt\+\sigma_2dW_{j}, \quad \cx_j(0)\=\eta_{j}, \\
   &dy_i \= -\left( A_1^Ty_i   \+ Q_1 \cx_i \+ S_1 \right) \+ z_idW_i,\quad y_i(T) \= 0,\\
   &dy_j \= -\left( A_2^Ty_j   \+ Q_2 \cx_j \+ S_2 \right) \+ z_jdW_j,\quad y_j(T) \= 0,\\
   &R_1\cu_i \+   B_1^Ty_i=0,\quad R_2\cu_j \+ B_2^Ty_j=0.
   \end{aligned}\right.
   \end{equation*}
\end{Proposition}

\subsection{Consistency condition}

This sub-step aims to synthesize \emph{all} generic behaviors $\{\check{x}_{i}(\check{u}_{i}(\textbf{m}, \textbf{P}^*))\}_{i \in \I_1}$ in LS1, $\{\check{x}_{j}(\check{u}_{j}(\textbf{m}, \widehat{\textbf{P}}^*))\}_{j \in \I_2}$ in LS2, to match aggregated $\textbf{m}$ across LS. Noting LS1 itself is homogenous (although LS is not), optimal auxiliary controls $\{\check{u}_{i}\}_{i \in \I_1}$  are thus symmetric and so is $\{\check{x}_{i}(\check{u}_{i})\}_{i \in \I_1}.$ Parallel holds for $\{\check{x}_{j}(\check{u}_{j})\}_{j \in \I_2}.$ Then, applying de Finetti's theorem to realized $\{\check{x}_{i}(\check{u}_{i}(\textbf{m}, \textbf{P}^*))\}_{i \in \I_1}$, $\{\check{x}_{j}(\check{u}_{j}(\textbf{m}, \widehat{\textbf{P}}^*))\}_{j \in \I_2}$, we yield a bilateral fixed-point CC equivalence (noticing that $\textbf{m}=(\m, \m_1, \m_2)$ ):
\begin{equation}\label{cc}
\text{\textbf{(CC)}}
\left\{\begin{aligned}
 \text{LS1:} \ \ &\m_1=\mathbb{E}(\bar{x}_{i}(\bar{u}_{i}(\m_1, \m_2, \textbf{P}^*))), \\
&\textbf{P}^*=\textbf{P}^*(\bar{x}_{i}(\bar{u}_{i}(\m_1, \m_2, \textbf{P}))); \\
 \text{LS2:} \ \ &\m_2=\mathbb{E}(\bar{x}_{j}(\bar{u}_{j}(\m_1, \m_2, \widehat{\textbf{P}}^*))),\\
&\widehat{\textbf{P}}^*=\widehat{\textbf{P}}^*(\bar{x}_{j}(\bar{u}_{j}(\m_1, \m_2, \widehat{\textbf{P}}^*))).
\end{aligned}
\right.
\end{equation}
So $\textbf{m}$ can be specified. The CC system is given by
\begin{equation}\label{ccsystem}
   \left\{\begin{aligned}
&d\bx = \left(\bA\bx \+ \bB\bu \+ \mathbf{F}\E\bx   \right)dt \+\boldsymbol{\sigma}_1
dW_i \+ \boldsymbol{\sigma}_2dW_j,\\%
&d\by = \left(\bar{\bA}\bx \+ \bar{\bB}\by  \+  \widetilde{\mathbf{F}}\E\bx \+ \widetilde{\mathbf{H}}\E\by  \right)dt \+ \bz_1dW_i  \+ \bz_2dW_j,\\
&\bR\bu \+ \widetilde{\mathbf{B}}^T\by  = 0,\ \ \bx(0)=\bx_0,\quad\by(T)=\mathbf{0},
\end{aligned}\right.
   \end{equation}
\n
where
\fontsize{4.9}{6}\selectfont
\begin{equation*}
\left\{\begin{aligned}
&\bx=(\cx_i,\cx_j),\ \ \by=(y_i,p^{(1)*},p^{(2)*},p^{(1)*}_i,p^{(2)*}_j,y_j,\hat{p}^{(2)*},\hat{p}^{(1)*}, \hat{p}^{(2)*}_j,\hat{p}^{(1)*}_i), \\
&\bx_0=(\xi_i,\eta_j),\ \bz_1=(z_i,0,0,q_i^{(1)*},0,0,0,0,0,{\hat{q}}_i^{(1)*}), \\
& \bz_2=(0,0,0,0,{q}_j^{(2)*},z_j,0,0,\hat{q}_j^{(2)*},0),\ \ \bu=(\cu_i,\cu_j),\\
&\bA =
\left(\begin{smallmatrix}
               A_1 & 0 \\
               0 & A_2
\end{smallmatrix}\right),
\bB =
\left(\begin{smallmatrix}
               B_1 & 0 \\
               0 & B_2
\end{smallmatrix}\right),
\mathbf{F} =
\left(\begin{smallmatrix}
               \pi_1F_1 & \pi_2F_1 \\
               \pi_1F_2 & \pi_2F_2
\end{smallmatrix}\right),
\mathbf{R} =
\left(\begin{smallmatrix}
               R_1 &  0 \\
                0 & R_2
\end{smallmatrix}\right),
\\
&\bar{\bA} =
\left(\begin{smallmatrix}
               -Q_1 & 0 \\
               0 & 0\\
               0 & 0\\
               -Q_1 & 0\\
               0 & -\alpha Q_2\\
               0 & -Q_2\\
               0 & 0\\
               0 & 0\\
               0 & -  Q_2\\
               -\beta Q_1 & 0\\
\end{smallmatrix}\right),
\tilde{\mathbf{B}} =
\left(\begin{smallmatrix}
               B_1 & 0 \\
               0 & 0\\
               0 & 0\\
              0 & 0\\
               0 & 0\\
               0 & B_2\\
               0 & 0\\
               0 & 0\\
               0 & 0\\
               0 & 0\\
\end{smallmatrix}\right),\boldsymbol{\sigma}_1=\left(\begin{smallmatrix}
  \sigma_1 \\
  0
\end{smallmatrix}\right),
\boldsymbol{\sigma}_2=\left(\begin{smallmatrix}
  0\\
  \sigma_2
\end{smallmatrix}\right),\\
&\bar{\bB} =
\left(\begin{smallmatrix}
               -A_1^T & - \pi_1 F_1^T & - \pi_2 F_2^T & 0 & 0& 0& 0& 0& 0& 0\\
               0 & - A_1^T \- \pi_1F_1^T &  -\pi_2F_2^T & 0 & 0& 0& 0& 0& 0& 0\\
               0 &  - \pi_1 F_1^T & - A_2^T \- \pi_2F_2^T & 0 & 0& 0& 0& 0& 0& 0\\
               0 & 0 & 0 & - A_1^T &  0 & 0 & 0 & 0 & 0 & 0\\
               0 & 0 & 0 & 0 & - A_2^T &  0 & 0 & 0 & 0 & 0\\
               0 & 0 & 0 & 0 & 0 & -A_2^T & - \pi_2 F_2^T & - \pi_1 F_1^T & 0& 0\\
               0 & 0 & 0 & 0 & 0 & 0 & - A_2^T \- \pi_2F_2^T &  -\pi_1 F_1^T & 0 &   0\\
               0 & 0& 0& 0& 0& 0& - \pi_2 F_2^T & - A_1^T \- \pi_1F_1^T &  0&  0\\
               0 & 0 & 0 & 0 &  0 & 0 & 0 & 0 & - A_2^T & 0\\
               0 & 0 & 0 & 0 & 0 &  0 & 0 & 0 & 0 & - A_1^T\\
\end{smallmatrix}\right),\\
&\tilde{\mathbf{F}}  =
\left(\begin{smallmatrix}
            \pi_1(Q_1\Gamma_1 \+ \Gamma_1^TQ_1 \-  \pi_1\Gamma_1^TQ_1\Gamma_1-\alpha\pi_2\Gamma_2^TQ_2\Gamma_2)   & \pi_2(Q_1\Gamma_1 \+ \alpha\Gamma_2^TQ_2 \-  \pi_1\Gamma_1^TQ_1\Gamma_1-\alpha\pi_2\Gamma_2^TQ_2\Gamma_2) \\
               \pi_1(Q_1\Gamma_1 \+ \Gamma_1^TQ_1 \-  \pi_1\Gamma_1^TQ_1\Gamma_1-\alpha\pi_2\Gamma_2^TQ_2\Gamma_2)   & \pi_2(Q_1\Gamma_1 \+ \alpha\Gamma_2^TQ_2 \-  \pi_1\Gamma_1^TQ_1\Gamma_1-\alpha\pi_2\Gamma_2^TQ_2\Gamma_2) \\
               \pi_1(\Gamma_1^TQ_1 \+\alpha Q_2\Gamma_2 \- \pi_1\Gamma_1^TQ_1\Gamma_1-\alpha\pi_2\Gamma_2^TQ_2\Gamma_2)   & \pi_2(\alpha Q_2\Gamma_2 \+ \alpha\Gamma_2^TQ_2 \-  \pi_1\Gamma_1^TQ_1\Gamma_1-\alpha\pi_2\Gamma_2^TQ_2\Gamma_2) \\
              0 & 0\\
               0 & 0\\
               \pi_1(Q_2\Gamma_2 \+\beta \Gamma_1^TQ_1 \-  \beta\pi_1\Gamma_1^TQ_1\Gamma_1-\pi_2\Gamma_2^TQ_2\Gamma_2)   & \pi_2(Q_2\Gamma_2 \+ \Gamma_2^TQ_2 \-  \beta\pi_1\Gamma_1^TQ_1\Gamma_1-\pi_2\Gamma_2^TQ_2\Gamma_2) \\
               \pi_1(Q_2\Gamma_2 \+\beta \Gamma_1^TQ_1 \-  \beta\pi_1\Gamma_1^TQ_1\Gamma_1-\pi_2\Gamma_2^TQ_2\Gamma_2)   & \pi_2(Q_2\Gamma_2 \+ \Gamma_2^TQ_2 \-  \beta\pi_1\Gamma_1^TQ_1\Gamma_1-\pi_2\Gamma_2^TQ_2\Gamma_2) \\
              \pi_1(\beta Q_1\Gamma_1 \+\beta \Gamma_1^TQ_1 \-  \pi_1\Gamma_1^TQ_1\Gamma_1-\pi_2\Gamma_2^TQ_2\Gamma_2)   & \pi_2(\beta Q_1\Gamma_1 \+ \Gamma_2^TQ_2 \-  \pi_1\Gamma_1^TQ_1\Gamma_1-\pi_2\Gamma_2^TQ_2\Gamma_2) \\
               0 & 0\\
               0 & 0\\
\end{smallmatrix}\right),\\
\end{aligned}\right.
\end{equation*}
\begin{equation*}
\begin{aligned}
&\tilde{\mathbf{H}}  =
\left(\begin{smallmatrix}
            0 & -\pi_1F_1^T  & -\pi_2F_2^T &  -\pi_1F_1^T &  -\pi_2F_2^T & 0& 0&0&0&0  \\
            0 & -\pi_1F_1^T  & -\pi_2F_2^T &  -A_1^T-\pi_1F_1^T &  -\pi_2F_2^T & 0& 0&0&0&0\\
            0 & -\pi_1F_1^T  & -A_2^T-\pi_2F_2^T &  -\pi_1F_1^T &  -\pi_2F_2^T & 0& 0&0&0&0\\
              0 & 0& 0 & 0 & 0& 0 & 0 & 0& 0 & 0\\
               0 & 0& 0 & 0 & 0& 0 & 0 & 0& 0 & 0\\
               0 & 0& 0& 0&0 & 0& -\pi_2F_2^T  & -\pi_1F_1^T &  -\pi_2F_2^T &  -\pi_1F_1^T  \\
           0 & 0& 0& 0&0 & 0& -A_2^T-\pi_2F_2^T  & -\pi_1F_1^T &  -\pi_2F_2^T &  -\pi_1F_1^T  \\
           0 & 0& 0& 0&0 & 0& -\pi_2F_2^T  & -A_1^T-\pi_1F_1^T &  -\pi_2F_2^T &  -\pi_1F_1^T  \\
              0 & 0& 0 & 0 & 0& 0 & 0 & 0& 0 & 0\\
               0 & 0& 0 & 0 & 0& 0 & 0 & 0& 0 & 0\\
\end{smallmatrix}\right).
\end{aligned}
\end{equation*}  \normalsize
and the mean field terms are determined by $\m_1 = \E\cx_i$, $\m_2 = \E\cx_j$.

\begin{Proposition}
  If Riccati equations
  \begin{equation}\label{RE_LY}
\left\{\begin{aligned}
   &\dot{\bP} \+  \bP\left(\bA \+ \mathbf{F}\right) \- \left(\widetilde{\mathbf{H}} \+ \bar{\bB}\right)\bP \- \bP\bB\bR^{-1}\tilde{\mathbf{B}}^T\bP \\
   &\qquad- \left(\bar{\bA}   \+ \widetilde{\mathbf{F}}\right) = 0,\\
   &  \dot{\bK} \+ \bK \left(\bA-\bB\bR^{-1}\tilde{\mathbf{B}}^T\bP\right) \- \left(\bar{\bB}  \+ \bP\bB\bR^{-1}\tilde{\mathbf{B}}^T\right)\bK \\
   &\qquad\+ \bK\bB\bR^{-1}\tilde{\mathbf{B}}^T  \bK \+ \bP\mathbf{F}  \-  \widetilde{\mathbf{F}}  \- \widetilde{\mathbf{H}}\bP  = 0,\\
&\bP(T) = 0,\quad \bK(T) = 0
\end{aligned}\right.
\end{equation}
  admit  solutions $\bP\!\in\! C([0,T];\mathbb{R}^{10n\times 2n})$, $\bK\!\in\! C([0,T];\mathbb{R}^{10n\times 2n})$, then Problem \ref{problem1} admits a feedback form decentralized control $\tilde{\bu} = \Theta_1\tilde{\bx} \+ \Theta_2\left(\E\tilde{\bx}\- \tilde{\bx}\right)$, where $\tilde{\bu}=(\tilde{u}_i,\tilde{u}_j)$, $\tilde{\bx}=(\tilde{x}_i,\tilde{x}_j)$, and
$$ \Theta_1= - \bR^{-1} \tilde{\mathbf{B}}^T\bP,\quad \Theta_2=  - \bR^{-1} \tilde{\mathbf{B}}^T\bK.$$
  The realized state $\tilde{\bx}=(\tilde{x}_i,\tilde{x}_j)$ satisfies the following dynamic
  \begin{equation*}
\left\{\begin{aligned}
&d\tilde{x}_i\= \left(A_1\tilde{x}_i\+B_1\tilde{u}_i \+ F_1 \tilde{x}^{(N)}\right )dt\+\sigma_1dW_{i}, \quad \tilde{x}_i(0)\=\xi_{i}, \\
   &d\tilde{x}_j\= \left(A_2\tilde{x}_j\+B_2\tilde{u}_j\+ F_2 \tilde{x}^{(N)} \right)dt\+\sigma_2dW_{j}, \quad \tilde{x}_j(0)\=\eta_{j}. \\
\end{aligned}\right.
\end{equation*}
\end{Proposition}
For the solvability of \eqref{RE_LY}, we have the following result
\begin{Proposition}
  Let $\Psi $, $\Phi $ be fundamental matrices $
  \Psi = \left(\begin{smallmatrix}
    \bA\+\mathbf{F}  & -\bB\bR^{-1}\tilde{\mathbf{B}}^T \\
    \bar{\bA}   \+ \widetilde{\mathbf{F}} & \widetilde{\mathbf{H}} \+ \bar{\bB}
  \end{smallmatrix}\right)
  $ and
  $\Phi =
  \left(\begin{smallmatrix}
    \bA-\bB\bR^{-1}\tilde{\mathbf{B}}^T\bP  & \bB\bR^{-1}\tilde{\mathbf{B}}^T \\
     -\left( \bP\mathbf{F}  \-  \widetilde{\mathbf{F}}  \- \widetilde{\mathbf{H}}\bP\right) & \bar{\bB}  \+ \bP\bB\bR^{-1}\tilde{\mathbf{B}}^T
  \end{smallmatrix}\right)
  $. If
  \begin{equation*}
  \left\{\begin{aligned}
  &\left[(0,I)e^{\Psi(T-t)}\binom{0}{I}\right]^{-1}\in L^1(0,T;\R^{10n\times 10n}),\\
  &\left[(0,I)e^{\Phi(T-t)}\binom{0}{I}\right]^{-1}\in L^1(0,T;\R^{10n\times 10n}),\\
  \end{aligned}\right.
  \end{equation*}
  then \eqref{RE_LY} admits unique solutions $\bP$, $\bK$ as
  \begin{equation}\label{explicit PK}
  \left\{\begin{aligned}
  &\bP(t) = -\left[(0,I)e^{\Psi(T-t)}\binom{0}{I}\right]^{-1} (0,I)e^{\Psi(T-t)}\binom{I}{0},\\
  &\bK(t) = -\left[(0,I)e^{\Phi(T-t)}\binom{0}{I}\right]^{-1} (0,I)e^{\Phi(T-t)}\binom{I}{0}.\\
  \end{aligned}\right.
  \end{equation}

  \end{Proposition}

\begin{Remark}
Actually, for different $\alpha$ and $\beta$, we may reach different kinds of mean field problems and find some interesting phenomena. Details can be referred to Appendix A.
\end{Remark}

\subsection{Distributed design}
Thus, we can conclude the following procedure of deriving the mean-field strategy:
\begin{description}
  \item[Step 1] By \eqref{explicit PK}, the solution $(\bP, \bK)$ of Riccati equations \eqref{RE_LY}  can be obtained.
  \item[Step 2] For any agent $\A_i$ (or $\A_j$) in $\I_1$ (or $\I_2$), besides its own information, it still need the information of another generic agent denoted by $\A_j$ (or $\A_i$) in its opposite team $\I_2$ (or $\I_1$). Then it can obtain its feedback form mean-field strategy by
      \begin{equation*}
      \begin{aligned}
      \binom{\tilde{u}_i}{\tilde{u}_j} = \Theta_1\binom{\tilde{x}_i}{\tilde{x}_j} \+ \Theta_2\left(\E\binom{\tilde{x}_i}{\tilde{x}_j}\- \binom{\tilde{x}_i}{\tilde{x}_j}\right),
      \end{aligned}
      \end{equation*}
      where   $ \Theta_1= - \bR^{-1} \tilde{\mathbf{B}}^T\bP$ and
  $\Theta_2=  - \bR^{-1} \tilde{\mathbf{B}}^T\bK$.
  \item[Step 3] The realized states $\tilde{x}_i$, $\tilde{x}_j$
      satisfy the following bilateral closed-loop system:\f
      \begin{equation*}
      \begin{aligned}
      d\binom{\tilde{x}_i}{\tilde{x}_j} = &
      \left\{\left[\begin{pmatrix}
        A_1 & 0 \\
        0 & A_2
      \end{pmatrix} \+
      \begin{pmatrix}
        B_1 & 0 \\
        0 & B_2
      \end{pmatrix}\left(\Theta_1 \- \Theta_2\right) \right]\binom{\tilde{x}_i}{\tilde{x}_j}\right.\\
       &\left.\+ \begin{pmatrix}
        B_1 & 0 \\
        0 & B_2
      \end{pmatrix}\Theta_2 \E\binom{\tilde{x}_i}{\tilde{x}_j} \+ \binom{F_1}{F_2}\tilde{x}^{(N)}\right\}dt\\
       &+ \binom{\sigma_1}{0}dW_i  \+ \binom{0}{\sigma_2}dW_j,\quad \binom{\tilde{x}_i}{\tilde{x}_j}(0) = \binom{\xi_i}{\eta_j}.
      \end{aligned}
      \end{equation*}
\end{description}
\begin{Remark}
It should be noticed that there are MF terms in \textbf{Step 2-3}. To break up this couple feather, we should take the expectation to the system, which becomes an ordinary differential equation (ODE). By solving the ODE, we derive $\E\binom{\tilde{x}_i}{\tilde{x}_j}$ which implies the MF term, and then we focus on the SDE without MF terms.
\end{Remark}

\section{Performance analysis} \label{sec 5}

We present \textbf{(Mix)} analysis in a three-step procedure by highlighting its novelty via pairwise comparisons to ``pure" MG and MT. All limits below are in $N \rightarrow  +\infty$ sense.

We continue to analyze the performance of \textbf{(Mix)} strategy $(\bar{\textbf{u}}_{1}, \bar{\textbf{u}}_{2})$ derived by Section \ref{sec 4}, in an asymptotic mixed-equilibrium-optima (\textbf{AMEO}) sense (it combines \emph{equilibrium} due to game, and (social) \emph{optima} due to team): \begin{equation}\label{amo}\text{\textbf{(AMEO )}}\left\{\begin{aligned}
&\exists\ \rho,\ \lambda>0, \text{s.t.} \ \ \text{LS1}:
\\ &\left|\frac{1}{N_{1}}\mathbf{J}_{\text{mix}}^{(1)}(\widehat{\textbf{u}}_{1},\bar{\textbf{u}}_{2})-\frac{1}{N_{1}}\mathbf{J}_{\text{mix}}^{(1)}
(\bar{\textbf{u}}_{1},\bar{\textbf{u}}_{2})\right|\\
&=O(N_{1}^{-\rho}+\epsilon_{N}^\lambda), \\
& \exists\ \rho,\ \lambda>0, \text{s.t.} \ \ \text{LS2}: \\ &\left|\frac{1}{N_{2}}\mathbf{J}_{\text{mix}}^{(1)}(\bar{\textbf{u}}_{1},\widehat{\textbf{u}}_{2})-\frac{1}{N_{2}}\mathbf{J}_{\text{mix}}^{(1)}
(\bar{\textbf{u}}_{1},\bar{\textbf{u}}_{2})\right|\\
&=O(N_{2}^{-\rho}+\epsilon_{N}^\lambda), \end{aligned}\right.\end{equation}
where $\epsilon_{N}=\sup\limits_{1\leq k\leq 2}\left|\pi_k^{(N)}-\pi_k\right|$.

AMEO poses an \emph{inside-outside-mixed} concept: social (team) optima inside, and Nash (game) equilibrium outside (LS1,\!\!\! \!\! LS2), along perturbed
$\widehat{\textbf{u}}_{1}:=(\cdots \hat{u}_{k} \cdots)|_{k \in \I_{1}},$ $\widehat{\textbf{u}}_{2}:=(\cdots \hat{u}_{k} \cdots)|_{k \in \I_{2}}.$ AMEO mainly includes three sub-steps  as below.
\begin{enumerate}
  \item Study convergence behavior of realized empirical LS1, LS2 averages to off-lined CC system (cf.\eqref{cc}) in an $L^{2}$-norm. A key point here is a refined bilateral forward-backward SDE (FBSDE) estimates on $(\bar{\textbf{u}}_{1},\bar{\textbf{u}}_{2}).$
  \item Estimate two upper bound(s) for candidate perturbation(s) $(\widehat{\textbf{u}}_{1}, \widehat{\textbf{u}}_{2})$ respectively subject to $\mathbf{J}_{\text{mix}}^{(1)}(\widehat{\textbf{u}}_{1},\bar{\textbf{u}}_{2}) \leq \mathbf{J}_{\text{mix}}^{(1)}
(\bar{\textbf{u}}_{1},\bar{\textbf{u}}_{2}), \mathbf{J}_{\text{mix}}^{(2)}(\bar{\textbf{u}}_{1},\widehat{\textbf{u}}_{2}) \leq \mathbf{J}_{\text{mix}}^{(2)}
(\bar{\textbf{u}}_{1},\bar{\textbf{u}}_{2}).$ Two keys: some tailor-made FBSDE stability estimates; $\mathbf{J}_{\text{mix}}^{(1)}$ convexity by underlying $(Q, R)$ positiveness in \eqref{JMG}.
  \item Formulate $\mathbf{J}_{\text{mix}}^{(1)}$ in \eqref{amo} as a \emph{quadratic} functional on $({\textbf{u}}_{1}, {\textbf{u}}_{2})$ by specifying its second-order-operator(s) for convexity, and first-order ones for gradient. Its Fr\'{e}chet derivatives to team-wise $(\textbf{u}_{1}, \textbf{u}_{2})$, and componentwise $\{u_{k}\}_{k \in \I}$ can thus be structured. Then, AMEO can be verified by combining $(\widehat{\textbf{u}}_{1}, \widehat{\textbf{u}}_{2})$ bounds in step (ii), a near-stationarity (asymptotic zero Fr\'{e}chet derivative) by step (i), and (uniform) \emph{convexity} of $\mathbf{J}_{\text{mix}}^{(1)}$.
\end{enumerate}

Now, we give the main result of this work, whose proof is based on some lemmas. Please refer to Appendix B-E for details.
\begin{Theorem}
  Under (H1)-(H4), the mean-field strategy $(\tilde{\bu}_1,\tilde{\bu}_2)$ satisfies the following asymptotic optimality
 \begin{equation}\label{nash asym}
  \left\{\begin{aligned}
  &\frac{1}{N_1}\left(\mathbf{J}_{\mathrm{mix}}^{(1)}( \tilde{\bu}_{1}, \tilde{\bu}_{2}) \- \inf_{{\bu}_{1}}\mathbf{J}_{\mathrm{mix}}^{(1)}( {\bu}_{1}, \tilde{\bu}_{2})\right) = O(N_1^{-\frac{1}{2}}+\epsilon_{N}),\\
  &\frac{1}{N_2}\left(\mathbf{J}_{\mathrm{mix}}^{(2)}( \tilde{\bu}_{1}, \tilde{\bu}_{2}) \- \inf_{{\bu}_{2}}\mathbf{J}_{\mathrm{mix}}^{(2)}( {\tilde\bu}_{1}, {\bu}_{2})\right) = O(N_2^{-\frac{1}{2}}+\epsilon_{N}).
  \end{aligned}\right.
  \end{equation}
  Hence, $(\tilde{\bu}_1,\tilde{\bu}_2)$ is an asymptotic Nash equilibrium for the game between $\I_1$ and $\I_2$.
\end{Theorem}

\emph{Proof}
  For $L^2$ bounded candidate $\hat{\bu}_1$, by letting $\d \bu_1 \= \tilde{\bu}_{1} \- \hat{\bu}_{1}$ we have
\begin{equation*}
\begin{aligned}
\mathbf{J}_{\mathrm{mix}}^{(1)}( \tilde{\bu}_{1}, \tilde{\bu}_{2}) \- \mathbf{J}_{\mathrm{mix}}^{(1)}( \hat{\bu}_{1}, \tilde{\bu}_{2}) =      2\langle \M_2\tilde{\bu}_1 \+ \M_1,\d\mathbf{u}_1 \rangle \+ o(\d\bu_1).
\end{aligned}
\end{equation*}
Here, $\langle \M_2\tilde{\bu}_1 + \M_1 , \cdot\rangle$ is the Fr\'{e}chet derivative of ${J}^{(N)}_{soc}$ on $\tilde{\mathbf{u}}_i$.
Due to the linearity, we also get
\begin{equation*}
 \begin{aligned}
 &\mathbf{J}_{\mathrm{mix}}^{(1)}( \tilde{\bu}_{1}, \tilde{\bu}_{2}) \- \mathbf{J}_{\mathrm{mix}}^{(1)}( \hat{\bu}_{1}, \tilde{\bu}_{2}) \\
  =\ &2\sum_{i=1}^{N_1}\langle \M_2\tilde{\mathbf{u}} + \M_1 , \delta\mathbf{u}_1^i\rangle + o(\delta \mathbf{u}), \\
 \end{aligned}
\end{equation*}
where $\delta \mathbf{u}_1^i := (0,\cdots,0, \tilde{u}_i \- \hat{u}_i, 0,\cdots,0)$. Based on the synthesis analysis in Section \ref{sec 3}, we derive
\begin{equation*}
\begin{aligned}
&\langle \M_2\tilde{\mathbf{u}} + \M_1 , \delta\mathbf{u}_1^i\rangle= \d \mathbf{J}_{\text{mix}}^{(1)}(\d u_{i}) \\
 =\ &\mathbb{E}\int_0^{T}\Big[ \left\langle\bar{\Theta}_{1} \+ \pi_1  F_1^T\E p_k^{(1)*}  \+ \pi_2  F_2^T\E p_k^{(2)*} \+  \pi_1F_1^Tp^{(1)*} \right.\\
\ &\left. \+ \pi_2F_2^Tp^{(2)*}, \d x_{i}\right\rangle+\left\langle\Theta_{2}  , \d u_{i}\right\rangle\Big] dt \+\varepsilon_1^i  \+ \varepsilon_2^i \+ \varepsilon_3^i \\
=\ &\varepsilon_1^i  \+ \varepsilon_2^i \+ \varepsilon_3^i.
\end{aligned}
\end{equation*}
By virtue of some FBSDE estimations, we obtain
\begin{equation*}
\begin{aligned}
0&\leq \mathbf{J}_{\mathrm{mix}}^{(1)}( \tilde{\bu}_{1}, \tilde{\bu}_{2}) \- \mathbf{J}_{\mathrm{mix}}^{(1)}( \hat{\bu}_{1}, \tilde{\bu}_{2}) \\
 &\leq  \sum_{i=1}^{N_1}\left|\varepsilon_1^i  \+ \varepsilon_2^i \+ \varepsilon_3^i\right| = N_1O(N_1^{-\frac{1}{2}}+\epsilon_{N}).
\end{aligned}
\end{equation*}
The first equation of \eqref{nash asym} is established. Applying similar argument, we obtain the second one.
\hfill$\Box$

\section{Conclusion} \label{conc}

This paper investigates a new class of ``mixed" mean-field analysis.
The $\mathbf{C}_4-$type \emph{coalition matrix} analysis is related to game-team problem. 
A novel bilateral person-by-person optimality is introduced and wellposedness of related CC system is investigated.
Game-team strategies are designed and a novel asymptotic mixed-equilibrium-optima is proposed. An interesting work for further study is to consider the $\mathbf{C}_5-$type game-team problem referred in this work.

\begin{frontmatter}
\title{Appendix}

\thanks[footnoteinfo]{This document supplies appendices of the paper ``Linear quadratic mean-field game-team analysis: a mixed coalition approach" by Huang, Qiu, Wang and Wu, submitted to Automatica.}

\end{frontmatter}

\textbf{Appendix A: Some cases for different $\alpha$ and $\beta$.}\\

(i)  Assume $\alpha \= \beta = 1$, then the cost functionals of Problem \ref{problem1} reduce to
\begin{equation*}
      \begin{aligned}
      & \mathbf{J}_{\text{mix}}^{(1)} = \mathbf{J}_{\text{mix}}^{(2)} = \mathcal{J}^{1,(N)}_{\text{soc}} \+   \mathcal{J}^{2,(N)}_{\text{soc}},\\
      \end{aligned}
      \end{equation*}
      which leads to a combination of traditional social optima problems.

       If we further assume $A_1 \= A_2 := A$, $B_1 \= B_2 := B$, $F_1 \= F_2 := F$, $Q_1 \= Q_2 := Q$, $R_1 \= R_2 := R$, $\Gamma_1 \= \Gamma_2 := \Gamma$, $\xi_i$ and $\eta_j$ ($1\leq i\leq N_1$, $1\leq j\leq N_2$) are i.i.d and denoted by $\zeta_k$ ($1\leq k\leq N$), in this case we derive $p^{(1)*} = {p}^{(2)*} = \hat{p}^{(1)*} = \hat{p}^{(2)*}:= \hat{p}^{ *}$. We also obtain that $(\x_i, \u_i, y_i, p^{(1)*}_i, p^{(2)*}_j)$, $(\x_j, \u_j, y_j, \hat{p}^{(1)*}_i, \hat{p}^{(2)*}_j)$ are homogeneous.

(ii)  Assume $\alpha \= \beta = -1$, then the cost functionals of Problem \ref{problem1} reduce to
\begin{equation*}
      \begin{aligned}
      & \mathbf{J}_{\text{mix}}^{(1)} = \mathcal{J}^{1,(N)}_{\text{soc}} -  \mathcal{J}^{2,(N)}_{\text{soc}},\quad  \mathbf{J}_{\text{mix}}^{(2)} = \mathcal{J}^{2,(N)}_{\text{soc}} -  \mathcal{J}^{1,(N)}_{\text{soc}},
      \end{aligned}
      \end{equation*}
      which leads to a combination of two social optima problems (\textbf{inside}) and a zero-sum problem (\textbf{outside}). This mix problem can be viewed as a development of two-person zero-sum game problem where two-person becomes two-team (e.g. \cite{Sun21}, etc ).
If we further apply conditions on the coefficients, one may simplify the CC system and details are omitted.

(iii)  Assume $\alpha \=1, \beta = 0$, or $\alpha \=0, \beta = 1$, then the cost functionals reduce to
\begin{equation*}
      \begin{aligned}
      & \mathbf{J}_{\text{mix}}^{(1)} = \mathcal{J}^{1,(N)}_{\text{soc}} \+ \mathcal{J}^{2,(N)}_{\text{soc}},\quad  \mathbf{J}_{\text{mix}}^{(2)} = \mathcal{J}^{2,(N)}_{\text{soc}},
      \end{aligned}
      \end{equation*}
      or
\begin{equation*}
      \begin{aligned}
      \mathbf{J}_{\text{mix}}^{(1)} = \mathcal{J}^{1,(N)}_{\text{soc}},\quad  \mathbf{J}_{\text{mix}}^{(2)} = \mathcal{J}^{2,(N)}_{\text{soc}}\+ \mathcal{J}^{1,(N)}_{\text{soc}},
      \end{aligned}
      \end{equation*}
which leads to two classes of social optima. It means one group focuses on a social optima in itself manner; while the other group would like to consider both intra-group and inter-group cooperations.

(iv)  Assume $\alpha \=-1, \beta = 0$, or $\alpha \=0, \beta =-1$, then the cost functionals reduce to
\begin{equation*}
      \begin{aligned}
      & \mathbf{J}_{\text{mix}}^{(1)} = \mathcal{J}^{1,(N)}_{\text{soc}} \- \mathcal{J}^{2,(N)}_{\text{soc}},\quad  \mathbf{J}_{\text{mix}}^{(2)} = \mathcal{J}^{2,(N)}_{\text{soc}},
      \end{aligned}
      \end{equation*}
      or
\begin{equation*}
      \begin{aligned}
      \mathbf{J}_{\text{mix}}^{(1)} = \mathcal{J}^{1,(N)}_{\text{soc}},\quad  \mathbf{J}_{\text{mix}}^{(2)} = \mathcal{J}^{2,(N)}_{\text{soc}}\- \mathcal{J}^{1,(N)}_{\text{soc}},
      \end{aligned}
      \end{equation*}
which means one group focuses on a social optima in itself manner; while the other group would like to collaborate within the group itself and compete with the former group.

It should be noticed that, as a special case if $N=2$, $\mathcal{J}^{1,(N)}_{\text{soc}}=\mathcal{J}_1,$ $\mathcal{J}^{2,(N)}_{\text{soc}}=\mathcal{J}_2,$ and a group referred above becomes an agent. Thus above problems turn out to be optimal control problem, two-person social optima, etc.
\\

\textbf{Appendix B: Estimation of mean field coupling.}\\
\begin{Lemma}\label{lemma 6.1}
 Under (H1)-(H4), it holds that
  \begin{equation*}
       \begin{aligned}
      & \Big\|\frac{1}{N_1}\sum_{k \in \I_1} \bar{x}_{k}-\m_1\Big\|^{2}_{L^{2}}=O(N_{1}^{-1}+\epsilon_{N}^2),\\
      & \Big\|\frac{1}{N_2}{\sum_{k \in \I_2} \bar{x}_{k}}-\m_2\Big\|^{2}_{L^{2}}=O(N_{2}^{-1}+\epsilon_{N}^2).
       \end{aligned}
       \end{equation*}
\end{Lemma}

\emph{Proof}
  The dynamic of the realized state average is
  \begin{equation}\label{tilde x N1 N2}
\left\{\begin{aligned}
&d  \tilde{x}^{(N_1)}\= \left(A_1\tilde{x}^{(N_1)}\+B_1\tilde{u}^{(N_1)} \+ F_1 \pi_1^{(N)}\tilde{x}^{(N_1)}\right.\\
&\left.\qquad\qquad + F_1 \pi_2^{(N)}\tilde{x}^{(N_2)}
  \right)dt\+ \frac{1}{N_1}\sum_{i \in \I_1} \sigma_1dW_{i}, \\
   &d\tilde{x}^{(N_2)}\= \left(A_2\tilde{x}^{(N_2)}\+B_2\tilde{u}^{(N_2)} \+ F_2 \pi_1^{(N)}\tilde{x}^{(N_1)}\right.\\
   &\left.\qquad\qquad+ F_2 \pi_2^{(N)} \tilde{x}^{(N_2)} \right)dt\+\frac{1}{N_2}\sum_{j \in \I_2} \sigma_2dW_{j},  \\
  & \tilde{x}^{(N_1)}(0)\=\xi_{i}, \quad\tilde{x}^{(N_2)}(0)\=\eta_{j},
\end{aligned}\right.
\end{equation}
where  $\tilde{x}^{(N_1)}:= \frac{1}{N_1}\sum_{k \in \I_1}  \tilde{x}_{k}$, $\tilde{x}^{(N_2)}:= \frac{1}{N_2}\sum_{k \in \I_2}  \tilde{x}_{k}$, $\tilde{u}^{(N_1)}:= \frac{1}{N_1}\sum_{k \in \I_1}  \tilde{u}_{k}$, $\tilde{u}^{(N_2)}:= \frac{1}{N_2}\sum_{k \in \I_2}  \tilde{u}_{k}$.
The dynamic of the mean field is
\begin{equation}\label{bar bm N1 N2}
\left\{\begin{aligned}
&d\m_1\= (A_1\m_1\+B_1\E\cu_i \+ F_1 \pi_1\m_1 \+ F_1 \pi_2\m_2 )dt, \\
&d\m_2\= (A_2\m_2\+B_2\E\cu_j\+ F_2 \pi_1 \m_1 \+ F_2 \pi_2\m_2 )dt,\\
&\m_1(0)\=\E\xi,  \quad \m_2(0)\=\E\eta.
\end{aligned}\right.
\end{equation}
By applying Cauchy inequality, BDG inequality to the difference of \eqref{tilde x N1 N2}-\eqref{bar bm N1 N2}, for some positive constant $K$ which is independent of $N_1$, $N_2$,  it holds that
\s\begin{equation*}
\left\{\begin{aligned}
&\E\sup_{0\leq s\leq t}\|\tilde{x}^{(N_1)}(s) - \m_1(s)\|^2 \\
\leq &K \int_{0}^{t}\E\sup_{0\leq r\leq s}\|\tilde{x}^{(N_1)}(r) - \m_1(r)\|^2ds \+ O(N_1^{-1}+\epsilon_{N}^2),\\
&\E\sup_{0\leq s\leq t}\|\tilde{x}^{(N_2)}(s) - \m_2(s)\|^2 \\
\leq &K \int_{0}^{t}\E\sup_{0\leq r\leq s}\|\tilde{x}^{(N_2)}(r) - \m_2(r)\|^2ds \+ O(N_2^{-1}+\epsilon_{N}^2).\\
\end{aligned}\right.
\end{equation*}\n
By applying Gr\"{o}nwall's inequality, we have
\begin{equation*}
\left\{\begin{aligned}
&\E\sup_{0\leq t\leq T}\|\tilde{x}^{(N_1)}(t) - \m_1(t)\|^2 = O(N_1^{-1}+\epsilon_{N}^2),\\
&\E\sup_{0\leq t\leq T}\|\tilde{x}^{(N_2)}(t) - \m_2(t)\|^2 = O(N_2^{-1}+\epsilon_{N}^2).
\end{aligned}\right.
\end{equation*}
\hfill$\Box$\\

\textbf{Appendix C: Perturbations estimation.}\\

\begin{Lemma}\label{lemma 6.2}
Under (H1)-(H4), for some positive constant $K$ which is independent on $N_1$, $N_2$, it holds that
\begin{equation*}
  \begin{aligned}
  \mathbf{J}_{\mathrm{mix}}^{(1)}(\tilde{\textbf{u}}_{1},\tilde{\textbf{u}}_{2})\leq NK,\quad \mathbf{J}_{\mathrm{mix}}^{(2)}(\tilde{\textbf{u}}_{1},\tilde{\textbf{u}}_{2})\leq NK.
  \end{aligned}
  \end{equation*}
\end{Lemma}

\emph{Proof}
  For $\mathbf{J}_{\mathrm{mix}}^{(1)}(\tilde{\textbf{u}}_{1},\tilde{\textbf{u}}_{2})$, we have the following decompositions
  \begin{equation*}
\left\{\begin{aligned}
&\mathcal{J}_{i}(\bu_1,\bu_2)\=\frac{1}{2}\mathbb{E}\int_0^{T}\Big[\|x_{i}-\Gamma_1 x^{(N)}\|_{Q_1}^{2}\+\|u_{i}\|^{2}_{R_1}\Big]dt,\\
&\mathcal{J}_{j}(\bu_1,\bu_2)\=\frac{1}{2}\mathbb{E}\int_0^{T}\Big[\|x_{j}-\Gamma_2 x^{(N)}\|_{Q_2}^{2}\+\|u_{j}\|^{2}_{R_2}\Big]dt,
\end{aligned}\right.
\end{equation*}
\s\begin{equation*}
  \begin{aligned}
  &\mathbf{J}_{\mathrm{mix}}^{(1)}(\tilde{\textbf{u}}_{1},\tilde{\textbf{u}}_{2}) = \mathcal{J}^{1,(N)}_{\text{soc}}(\tilde{\textbf{u}}_{1},\tilde{\textbf{u}}_{2}) \+ \alpha \mathcal{J}^{2,(N)}_{\text{soc}}(\tilde{\textbf{u}}_{1},\tilde{\textbf{u}}_{2})\\
=\ & \frac{1}{2}\sum_{i=1}^{N_1}\mathbb{E}\int_0^{T}\Big[\|\tx_{i}-\Gamma_1 \tx^{(N)}\|_{Q_1}^{2}\+\|\tu_{i}\|^{2}_{R_1}\Big]dt\\
 &\+ \frac{\alpha}{2}\sum_{j=1}^{N_2}\mathbb{E}\int_0^{T}\Big[\|\tx_{j}-\Gamma_2 \tx^{(N)}\|_{Q_2}^{2}\+\|\tu_{j}\|^{2}_{R_2}\Big]dt\\
   \end{aligned}
\end{equation*}
\begin{equation*}
  \begin{aligned}
  \leq &\ \frac{K}{2}\left(\sum_{i=1}^{N_1}\mathbb{E}\int_0^{T}\Big[\|\tx_{i}- \m\|^{2} \+ \|\m \- \tx^{(N)}\|^{2} \+\|\tu_{i}\|^{2}\Big]dt\right. \\
  &\quad\left.\+\   \alpha\sum_{j=1}^{N_2}\mathbb{E}\int_0^{T}\Big[\|\tx_{j}- \m\|^{2} \+ \|\m \- \tx^{(N)}\|^{2} \+\|\tu_{j}\|^{2}\Big]dt\right).
  \end{aligned}
\end{equation*}\n
  By applying Cauchy inequality, BDG inequality and Gr\"{o}nwall's inequality, for some positive constant $K$ independent on $N_1$, $N_2$ we have
  \begin{equation*}
  \left\{\begin{aligned}
  &\sup_{i\in\I_1}\E\sup_{0\leq t \leq T} \|\tx_{i}(t)- \m(t)\|^2\leq K,\\
  &\sup_{j\in\I_2}\E\sup_{0\leq t \leq T} \|\tx_{j}(t)- \m(t)\|^2\leq K.\\
  \end{aligned}\right.
  \end{equation*}
  Combined with Lemma \ref{lemma 6.1}, it holds that
  \begin{equation*}
  \begin{aligned}
  \mathbf{J}_{\mathrm{mix}}^{(1)}(\tilde{\textbf{u}}_{1},\tilde{\textbf{u}}_{2})\leq NK.
  \end{aligned}
  \end{equation*}
  Similar argument can be applied to $\mathbf{J}_{\mathrm{mix}}^{(2)}$ and we complete the proof.
\hfill$\Box$\\

\textbf{Appendix D: Perturbations estimation.}\\
\begin{Lemma}
Under (H1)-(H4), since we are studying the asymptotic optimality of $(\tilde{\textbf{u}}_{1},\tilde{\textbf{u}}_{2})$, it is sufficient only to consider those admissible controls $(\hat{\textbf{u}}_{1},\hat{\textbf{u}}_{2})$ performing better than $(\tilde{\textbf{u}}_{1},\tilde{\textbf{u}}_{2})$. Specifically, $(\hat{\textbf{u}}_{1},\hat{\textbf{u}}_{2})$ satisfies
  \begin{equation}\label{candidate}
 \begin{aligned}
 \sum_{i=1}^{N_1}\mathbb{E}\int_{0}^{T}\|\hat{u}_i\|^2dt\leq NK, \quad \sum_{j=1}^{N_2}\mathbb{E}\int_{0}^{T}\|\hat{u}_j\|^2dt \leq NK. \\
 \end{aligned}
 \end{equation}
\end{Lemma}

\emph{Proof}
 Since $\hat{\textbf{u}}_{1}$ performs better then $\tilde{\textbf{u}}_{1}$, we have
  \begin{equation*}
  \begin{aligned}
  \sum_{i=1}^{N_1}\mathbb{E}\int_{0}^{T}\|\hat{u}_i\|^2dt &\leq \mathcal{J}^{1,(N)}_{\text{soc}}(\hat{\textbf{u}}_{1},\tilde{\textbf{u}}_{2}) \+ \alpha \mathcal{J}^{2,(N)}_{\text{soc}}(\hat{\textbf{u}}_{1},\tilde{\textbf{u}}_{2}) \\
  &=  \mathbf{J}_{\mathrm{mix}}^{(1)}(\hat{\textbf{u}}_{1},\tilde{\textbf{u}}_{2})\leq \mathbf{J}_{\mathrm{mix}}^{(1)}(\tilde{\textbf{u}}_{1},\tilde{\textbf{u}}_{2}) \leq NK.
  \end{aligned}
  \end{equation*}
   Similar argument can be applied to $\hat{\textbf{u}}_{2}$ and the proof is complete.
\hfill$\Box$\\

\textbf{Appendix E: Quadratic representation and Fr\'{e}chet derivatives.}\\

For the sake of notation simplicity, for a matrix $M$ and positive number $n_1$ and $n_2$, we introduce the following notations
\s\begin{equation*}
\left\{\begin{aligned}
&
\diag_{n_1}(M): =
\begin{smallmatrix}
  1 \\
  \vdots \\
  n_1
\end{smallmatrix}
\left(\begin{smallmatrix}
                  M &  & 0  \\
                    & \ddots\\
                   0 &  & M
                \end{smallmatrix}\right)\text{denotes an $n_1\times n_1$ size block}\\
               & \qquad\qquad\qquad\qquad\qquad\qquad\qquad \text{diagonal matrix}.\\
&\mathds{1}_{n_1\times n_2 }(M): =
\begin{smallmatrix}
  1 \\
  \vdots \\
  n_1
\end{smallmatrix}
\left(\begin{smallmatrix}
                  M & \dots & M  \\
                  \vdots & \ddots &\vdots\\
                   M & \dots & M
                \end{smallmatrix}\right)\text{denotes an $n_1\times n_2$ size block}\\
&\hspace{2.8cm}\begin{smallmatrix}
  1 \  \cdots\ \  n_2
\end{smallmatrix}\\
& \qquad\qquad\qquad\qquad\qquad\qquad\qquad \text{matrix generated by $M$}.\\
\end{aligned}\right.
\end{equation*}\n

Rewrite the problem as the following high-dimensional type
\s\begin{equation*}
\left\{\begin{aligned}
&d\mathbb{X} = \left(\mathbb{A} \mathbb{X} \+ \bbB_1 \bu_1 \+ \bbB_2 \bu_2  \right) dt \+ {\sum_{i=1}^{N_1}} \boldsymbol{\sigma}_i^1 dW_i\+  {\sum_{j=1}^{N_2}}\boldsymbol{\sigma}_j^2 dW_j,\\[-2mm]
&\mathbf{J}_{\mathrm{mix}}^{(1)}( \bu_{1}, \bu_{2}) \= \E {\int_{0}^{T}} \Big[\left\langle \bbQ_1 \bx,\bx\right\rangle \+  \left\langle \bbR_1\bu_1,\bu_1\right\rangle \+  \alpha\left\langle \bbR_2\bu_2,\bu_2\right\rangle \Big] dt,\\
&\mathbf{J}_{\mathrm{mix}}^{(2)}( \bu_{1}, \bu_{2}) \= \E {\int_{0}^{T} } \Big[ \left\langle {\bbQ}_2\bx,\bx\right\rangle \+  \left\langle  {\bbR}_2\bu_2,\bu_2\right\rangle \+ \beta \left\langle  {\bbR}_1\bu_1,\bu_1\right\rangle\Big]dt,
\end{aligned}\right.
\end{equation*}\n
where
\tiny\begin{equation*}
  \left\{\begin{aligned}
     & \mathbb{X}=(\bx_1,\bx_2), \bu_1 = (u_{\theta_1},\cdots,u_{\theta_{N_1}}), \bu_2 = (u_{\vartheta_1},\cdots,u_{\vartheta_{N_2}}),\\
    & \mathbb{A}=
    \left(\begin{smallmatrix}
            \diag_{N_1}(A_1)\+\mathds{1}_{N_1\times N_1}(\frac{F_1}{N}) & \mathds{1}_{N_1\times N_2}(\frac{F_1}{N}) \\
            \mathds{1}_{N_2\times N_1}(\frac{F_2}{N}) & \diag_{N_2\times N_2}(A_2)
          \end{smallmatrix}\right),
    \boldsymbol{\sigma}_i^1=\begin{smallmatrix}
      1 \\
      \vdots \\
      i \\
      \vdots \\
      N_1\\
            \vdots \\
      N \\
    \end{smallmatrix}
    \left(\begin{smallmatrix}
      0 \\
      \vdots\\
      \sigma_1\\
      \vdots \\
         0 \\
      \vdots\\
      0\\
    \end{smallmatrix}\right),
    \boldsymbol{\sigma}_j^2=\begin{smallmatrix}
      1 \\
      \vdots \\
        N_1\\
        \vdots \\
        j \\
      \vdots \\
      N \\
    \end{smallmatrix}
    \left(\begin{smallmatrix}
      0 \\
      \vdots\\
      0\\
      \vdots \\
        \sigma_2 \\
      \vdots\\
      0\\
    \end{smallmatrix}\right),\\
  &    \bbB_1=
    \left(\begin{smallmatrix}
      \diag_{N_1}(B_1) \\
        0
    \end{smallmatrix}\right),
    \bbB_2=
    \left(\begin{smallmatrix}
    0\\
      \diag_{N_2}(B_2) \\
    \end{smallmatrix}\right),
   \mathbb{R}_1 = \diag_{N_1}(R_1),
    \mathbb{R}_2 = \diag_{N_2}(R_2),\\
    & \mathbb{Q}_1 = \begin{pmatrix}
                       \mathbb{Q}_1^{(1,1)} & \mathbb{Q}_1^{(1,2)} \\
                       \mathbb{Q}_1^{(2,1)} & \mathbb{Q}_1^{(2,2)}
                     \end{pmatrix}, \mathbb{Q}_2 = \begin{pmatrix}
                       \mathbb{Q}_2^{(1,1)} & \mathbb{Q}_2^{(1,2)} \\
                       \mathbb{Q}_2^{(2,1)} & \mathbb{Q}_2^{(2,2)}
                     \end{pmatrix},\\
     & \mathbb{Q}_1^{(1,1)} = \diag_{N_1}(Q_1) \+ \frac{1}{N}\mathds{1}_{N_1\times N_1 }\left(\pi_1^{(N)}\Gamma_1^TQ_1\Gamma_1 \+ \Gamma_1^TQ_1 \+ Q\Gamma_1 \+ \alpha\frac{\pi_2^{(N)}}{N}\Gamma_2^TQ_2\Gamma_2\right), \\
    &\mathbb{Q}_1^{(1,2)} =  \frac{1}{N}\mathds{1}_{N_1\times N_2 } \left(\pi_1^{(N)}\Gamma_1^TQ_1\Gamma_1 \+ \alpha\Gamma_2^TQ_2 \+ Q\Gamma_1 \+ \alpha\frac{\pi_2^{(N)}}{N}\Gamma_2^TQ_2\Gamma_2\right),  \\
     &\mathbb{Q}_1^{(2,1)} =  \frac{1}{N}\mathds{1}_{N_2\times N_1 } \left(\pi_1^{(N)}\Gamma_1^TQ_1\Gamma_1 \+ \alpha Q_2\Gamma_2  \+ \Gamma_1^TQ \+ \alpha\frac{\pi_2^{(N)}}{N}\Gamma_2^TQ_2\Gamma_2\right), \\
     & \mathbb{Q}_1^{(2,2)} = \alpha\diag_{N_2}(Q_2) \+ \frac{1}{N}\mathds{1}_{N_2\times N_2 } \left(\pi_1^{(N)}\Gamma_1^TQ_1\Gamma_1 \+ \alpha\Gamma_1^TQ_2 \+ \alpha Q\Gamma_2 \+ \alpha\frac{\pi_2^{(N)}}{N}\Gamma_2^TQ_2\Gamma_2 \right),\\
     & \mathbb{Q}_2^{(1,1)} = \beta\diag_{N_1}(Q_1) \+ \frac{1}{N}\mathds{1}_{N_1\times N_1 }\left(\beta\pi_1^{(N)}\Gamma_1^TQ_1\Gamma_1 \+ \beta\Gamma_1^TQ_1 \+ \beta Q\Gamma_1 \+  \frac{\pi_2^{(N)}}{N}\Gamma_2^TQ_2\Gamma_2\right),\\
     &\mathbb{Q}_2^{(1,2)} =  \frac{1}{N}\mathds{1}_{N_1\times N_2 } \left(\beta\pi_1^{(N)}\Gamma_1^TQ_1\Gamma_1 \+  \Gamma_2^TQ_2 \+ \beta Q\Gamma_1 \+  \frac{\pi_2^{(N)}}{N}\Gamma_2^TQ_2\Gamma_2\right),  \\
     &\mathbb{Q}_2^{(2,1)} =   \frac{1}{N}\mathds{1}_{N_2\times N_1 }  \frac{1}{N}\left(\beta\pi_1^{(N)}\Gamma_1^TQ_1\Gamma_1 \+   Q_2\Gamma_2  \+ \beta\Gamma_1^TQ \+  \frac{\pi_2^{(N)}}{N}\Gamma_2^TQ_2\Gamma_2\right), \\
     & \mathbb{Q}_2^{(2,2)} =  \diag_{N_2}(Q_2)  \+ \frac{1}{N}\mathds{1}_{N_2\times N_2 }\left(\beta\pi_1^{(N)}\Gamma_1^TQ_1\Gamma_1 \+  \Gamma_1^TQ_2 \+ Q\Gamma_2 \+  \frac{\pi_2^{(N)}}{N}\Gamma_2^TQ_2\Gamma_2\right).
  \end{aligned}\right.
\end{equation*}\n
Because of the game structure between $\I_1$ and $\I_2$, for agents in $\I_1$,  they are facing with the following social optima problem
\begin{equation*}
\left\{\begin{aligned}
&d\mathbb{X} = \left(\bbA \mathbb{X} \+ \bbB_1 \bu_1 \+ \bbB_2 \tilde{\bu}_2  \right) dt \+ \sum_{i=1}^{N_1} \boldsymbol{\sigma}_i^1 dW_i\+ \sum_{j=1}^{N_2}\boldsymbol{\sigma}_j^2 dW_j,\\
&\min_{\bu_1}\ \mathbf{J}_{\mathrm{mix}}^{(1)}( \bu_{1}, \tilde{\bu}_{2}) \= \E\int_{0}^{T} \Big[\left\langle \bbQ_1 \mathbb{X},\mathbb{X}\right\rangle \+  \left\langle \bbR_1\bu_1,\bu_1\right\rangle \\
&\hspace{4cm}\+  \left\langle \alpha\bbR_2\widetilde{\bu}_2,\widetilde{\bu}_2\right\rangle \Big] dt.\\
\end{aligned}\right.
\end{equation*}
Since some bounded linear operators $(\M_2, \M_1, \M_0)$ are only dependent on the coefficients and $\tilde{\bu}_{2}$, we can rewrite $\mathbf{J}_{\mathrm{mix}}^{(1)}$ as the following quadratic form
\begin{equation*}
\begin{aligned}
\mathbf{J}_{\mathrm{mix}}^{(1)}( \bu_{1}, \tilde{\bu}_{2}) =  \langle \M_2\mathbf{u}_1,\mathbf{u}_1\rangle + 2\langle \M_1,\mathbf{u}_1 \rangle + \M_0.
\end{aligned}
\end{equation*}


\begin{thebibliography}{99}


\bibitem{ad2010} D. Andersson and B. Djehiche (2010). A maximum principle for SDEs of mean-field type. Applied Mathematics and Optimization, 63(3), 341-356.


\bibitem{bsyy2016} A. Bensoussan, K. Sung, S. Yam, and S. Yung (2016). Linear-quadratic mean field games. Journal of Optimization Theory and Applications, 169, 496-529.

%



\bibitem{blp2014} R. Buckdhan, J. Li, and S. Peng (2014). Nonlinear stochastic differential games involving a major player and a large number of collectively acting minor agents. SIAM Journal on Control and Optimization, 52, 451-492.

\bibitem{BC98} P. J. Buckley and M. C. Casson (1998). Analyzing foreign market entry strategies: extending the internalization approach. Journal of International Business Studies, 29(3), 539-561.



\bibitem{c2010} P. Cardaliaguet (2010). Notes on Mean Field Games, Technical report.

\bibitem{cd2013}
R. Carmona and F. Delarue (2013). Probabilistic analysis of mean-field games. SIAM Journal on Control and Optimization, 51, 2705-2734.






\bibitem{DHM20} X. Du, P. He and J. R. Martins(2020). A B-spline-based generative adversarial network model for fast interactive airfoil aerodynamic optimization. In AIAA Scitech 2020 Forum (p. 2128).







\bibitem{gs2014} D. Gomes and J. Saude (2014). Mean field games models--a brief survey. Dyn. Games Appl., 4, 110-154.



\bibitem{HHN18} Y. Hu, J. Huang and T. Nie (2018). Linear-quadratic-Gaussian mixed mean-field games with heterogeneous input constraints. SIAM Journal on Control and Optimization, 56, 2835-2877.





\bibitem{hcm2007a} M. Huang, P. E. Caines, and R. P. Malham\'{e} (2007). Large-population cost-coupled LQG problems with non-uniform agents: individual-mass behavior and decentralized $\varepsilon$-Nash equilibria. IEEE Transactions on Automatic Control, 52, 1560-1571.


\bibitem{hcm2012} M. Huang, P. E. Caines, and R. P. Malham\'{e} (2012). Social optima in mean-field LQG control: centralized and decentralized strategies. IEEE Transactions on Automatic Control, 57, 1736-1751.




\bibitem{ll2007} J. Lasry and P. Lions (2007). Mean field games. Japanese Journal of Mathematics, 2, 229-260.





\bibitem{MTMW09} T. Michalak, J. Tyrowicz, P. McBurney and M. Wooldridge (2009). Exogenous coalition formation in the e-marketplace based on geographical proximity. Electronic Commerce Research and Applications, 8(4), 203-223.

\bibitem{mb2018}
J.~Moon and T.~Ba\c{s}ar (2018). Linear quadratic mean field stackelberg differential games. Automatica, 97, 200-213.

\bibitem{NNY2020} S. Nguyen, D. Nguyen and G. Yin (2020). A stochastic maximum principle for switching diffusions using conditional mean-fields with applications to control problems. ESAIM: Control, Optimisation and Calculus of Variations, 26, 1-26.

\bibitem{ncmh2013} M. Nourian, P. E. Caines, R. P. Malham\'{e}, and M. Huang (2013). Nash, social and centralized solutions to consensus problems via mean-field control theory. IEEE Transactions on Automatic Control, 58, 639-653.

\bibitem{nm2018} G. Nuno and B. Moll (2018). Social optima in economics with heterogeneous agents. Review of Economic Dynamics, 28, 150-180.



\bibitem{p2016} H. Pham (2016). Linear quadratic optimal control of conditional McKean-Vlasov equation with random coefficients and applications. Probab. Uncertain. Quantita. Risk, 1, 1-26.



\bibitem{Sun21} J. Sun (2021). Two-Person Zero-Sum Stochastic Linear-Quadratic Differential Games. SIAM Journal on Control and Optimization, 59(3), 1804-1829.


\bibitem{STY08} N. Sun, W. Trockel and Z. Yang (2008). Competitive outcomes and endogenous coalition formation in an n-person game. Journal of Mathematical Economics, 44, 853-860.




\bibitem{tzb2014} H. Tembine, Q. Zhu, and T. Ba\c{s}ar (2014). Risk-sensitive mean-field games. IEEE Transactions on Automatic Control, 59, 835-850.


\bibitem{wz2017}  B. Wang, J. Zhang (2017). Social optima in mean-field linear-quadratic-Gaussian models with Markov jump parameters. SIAM Journal on Control and Optimization, 55, 429-456.





\bibitem{Z94} L. Zhou (1994). A new bargaining set of an n-person game and endogenous coalition formation. Games and Economic Behavior, 6(3), 512-526.



\end{thebibliography}
\end{document}